\newif\ifprint 
\newif\ifdraft
\newif\ifonly
\newif\ifsuppress
\newif\ifhyper
\newif\ifstandalone
\newif\ifshowprivate 
\newif\ifbells 

\printfalse
\onlyfalse
\suppresstrue
\hypertrue
\showprivatefalse
\bellstrue
\standalonetrue
\draftfalse

\def\path{}
\newcommand\relinput[2][]{{\ifx{}#1\else\edef\path{\path#1/}\fi\input{\path#2}}}
\newcommand\relinclude[2][]{{\ifx{}#1\else\edef\path{\path#1/}\fi\include{\path#2}}}

\ifsuppress
	\gdef\wlog#1{}
	\global\let\GenericInfo\relax
	\DeclareRobustCommand{\GenericInfo}[2]{}
\fi

\ifstandalone
	\documentclass{amsip-p-l}
	
	\def\cAMSIP{}
	\copyrightinfo{}{}
\else
	\documentclass{amsip-m-l}
\fi

\draftfalse
\ifdraft
	\calclayout
\fi

\usepackage{amssymb,latexsym,amsmath,amsthm}
	
\usepackage{xspace}
\usepackage{color}
\usepackage{diagrams} 
\usepackage[pdftex]{graphicx}

\ifstandalone
	\def\pdfsyncstop{}
	\def\pdfsyncstart{}
\else
\def\pdfsyncstop{}
\def\pdfsyncstart{}
\fi

\def\artfolder{art/}

\ifhyper
\ifprint
	\definecolor{linkred}{rgb}{0.2,0.2,0.2} 
	\definecolor{linkblue}{rgb}{0.2,0.2,0.2} 
\else
	\definecolor{linkred}{rgb}{0.7,0.2,0.2}
	\definecolor{linkblue}{rgb}{0,0.2,0.6}
\fi
\usepackage[
	hyperindex,
	pagebackref,
	pdftex,
	pdfusetitle,
	breaklinks=true,
	bookmarks=false,
	plainpages=false,
	colorlinks,
	linkcolor=linkblue,
	citecolor=linkblue,
	urlcolor=linkred,
	pdftoolbar=false,
	pdfpagelabels,
	pdfstartview=Fit
]
{hyperref}
\usepackage{memhfixc}
\usepackage[backrefs,msc-links,nobysame]{amsrefs}
\usepackage[low-sup]{subdepth}

\theoremstyle{definition}

\theoremstyle{remark}

\ifstandalone
\else
\numberwithin{section}{chapter}
\numberwithin{equation}{chapter}
\fi

\def\setdrawbox#1#2#3{
\pdfximage{#2} 
\setbox0=\hbox{\pdfrefximage\pdflastximage} 
\drawx=#1\wd0
\ifdim\drawx>\hsize\drawx=\hsize\fi
\pdfximage width \drawx {#2}
\setbox\drawing=\vbox{\offinterlineskip\pdfrefximage\pdflastximage\kern 0pt}
\drawx=\wd\drawing
\drawy=\ht\drawing
\ngap=0pt \sgap=0pt \wgap=0pt \egap=0pt
\setbox0=\vbox{\offinterlineskip \box\drawing \ifgridlines\drawgrid\drawx\drawy\fi #3}%
\setbox\drawing=\vbox{\kern\ngap\hbox{\kern\wgap\box0\kern\egap}\kern\sgap}
}
\def\obsoletedrawbox#1#2#3{\vbox{
  \setbox\drawing=\vbox{\offinterlineskip\epsfbox{#2.eps}\kern 0pt}
  \drawbp=\epsfurx
  \advance\drawbp by-\epsfllx\relax
  \multiply\drawbp by #1
  \divide\drawbp by 100
  \drawx=\drawbp truebp
  \ifdim\drawx>\hsize\drawx=\hsize\fi
  \epsfxsize=\drawx
  \setbox\drawing=\vbox{\offinterlineskip\epsfbox{#2.eps}\kern 0pt}
  \drawx=\wd\drawing
  \drawy=\ht\drawing
  \ngap=0pt \sgap=0pt \wgap=0pt \egap=0pt
  \setbox0=\vbox{\offinterlineskip
    \box\drawing \ifgridlines\drawgrid\drawx\drawy\fi #3}
  \kern\ngap\hbox{\kern\wgap\box0\kern\egap}\kern\sgap}}

\def\drawsmash#1{\relax
\ifmmode%
\typeout{Warning:Macro drawsmash used in math mode}%
\fi%
\setbox0=\hbox{#1}\ht0=0pt\dp0=0pt\box0}

\newcount\figcount
\newbox\drawing
\newcount\drawbp
\newdimen\drawx
\newdimen\drawy
\newdimen\ngap
\newdimen\sgap
\newdimen\wgap
\newdimen\egap

\newif\ifgridlines
\newbox\figtbox
\newbox\figgbox
\newdimen\figtx
\newdimen\figty

\newdimen\bwd
\def\hhline#1{\vbox{\drawsmash{\hbox to #1{\leaders\hrule height \bwd\hfil}}}}
\def\vvline#1{\hbox to 0pt{%
  \hss\vbox to #1{\leaders\vrule width \bwd\vfil}\hss}}

\def\clap#1{\hbox to 0pt{\hss#1\hss}}
\def\vclap#1{\vbox to 0pt{\offinterlineskip\vss#1\vss}}

\def\hstutter#1#2{\hbox{%
  \setbox0=\hbox{#1}%
  \hbox to #2\wd0{\leaders\box0\hfil}}}

\def\vstutter#1#2{\vbox{
  \setbox0=\vbox{\offinterlineskip #1}
  \dp0=0pt
  \vbox to #2\ht0{\leaders\box0\vfil}}}

\def\crosshairs#1#2{
  \dimen1=.002\drawx
  \dimen2=.002\drawy
  \ifdim\dimen1<\dimen2\dimen3\dimen1\else\dimen3\dimen2\fi
  \setbox1=\vclap{\vvline{2\dimen3}}
  \setbox2=\clap{\hhline{2\dimen3}}
  \setbox3=\hstutter{\kern\dimen1\box1}{4}
  \setbox4=\vstutter{\kern\dimen2\box2}{4}
  \setbox1=\vclap{\vvline{4\dimen3}}
  \setbox2=\clap{\hhline{4\dimen3}}
  \setbox5=\clap{\copy1\hstutter{\box3\kern\dimen1\box1}{6}}
  \setbox6=\vclap{\copy2\vstutter{\box4\kern\dimen2\box2}{6}}
  \setbox1=\vbox{\offinterlineskip\box5\box6}
  \drawsmash{\vbox to #2{\hbox to #1{\hss\box1}\vss}}}

\def\boxgrid#1{\rlap{\vbox{\offinterlineskip
  \setbox0=\hhline{\wd#1}
  \setbox1=\vvline{\ht#1}
  \drawsmash{\vbox to \ht#1{\offinterlineskip\copy0\vfil\box0}}
  \drawsmash{\vbox{\hbox to \wd#1{\copy1\hfil\box1}}}}}}

\def\drawgrid#1#2{\vbox{\offinterlineskip
  \dimen0=\drawx
  \dimen1=\drawy
  \divide\dimen0 by 10
  \divide\dimen1 by 10
  \setbox0=\hhline\drawx
  \setbox1=\vvline\drawy
  \drawsmash{\vbox{\offinterlineskip
    \copy0\vstutter{\kern\dimen1\box0}{10}}}
  \drawsmash{\hbox{\copy1\hstutter{\kern\dimen0\box1}{10}}}}}

\def\figtext#1#2#3#4#5{
  \setbox\figtbox=\hbox{#5}
  \dp\figtbox=0pt
  \figtx=-#3\wd\figtbox \figty=-#4\ht\figtbox
  \advance\figtx by #1\drawx \advance\figty by #2\drawy
  \dimen0=\figtx \advance\dimen0 by\wd\figtbox \advance\dimen0 by-\drawx
  \ifdim\dimen0>\egap\global\egap=\dimen0\fi
  \dimen0=\figty \advance\dimen0 by\ht\figtbox \advance\dimen0 by-\drawy
  \ifdim\dimen0>\ngap\global\ngap=\dimen0\fi
  \dimen0=-\figtx
  \ifdim\dimen0>\wgap\global\wgap=\dimen0\fi
  \dimen0=-\figty
  \ifdim\dimen0>\sgap\global\sgap=\dimen0\fi
  \drawsmash{\rlap{\vbox{\offinterlineskip
    \hbox{\hbox to \figtx{}\ifgridlines\boxgrid\figtbox\fi\box\figtbox}
    \vbox to \figty{}
    \ifgridlines\crosshairs{#1\drawx}{#2\drawy}\fi
    \kern 0pt}}}}

\def\swtext#1#2#3{\figtext{#1}{#2}00{#3}}
\def\setext#1#2#3{\figtext{#1}{#2}10{#3}}

\def\wtext#1#2#3{\figtext{#1}{#2}0{.5}{#3}}
\def\etext#1#2#3{\figtext{#1}{#2}1{.5}{#3}}
\def\ntext#1#2#3{\figtext{#1}{#2}{.5}1{#3}}
\def\stext#1#2#3{\figtext{#1}{#2}{.5}0{#3}}

\def\hpad#1#2#3{\hbox{\kern #1\hbox{#3}\kern #2}}
\def\vpad#1#2#3{\setbox0=\hbox{#3}\dp0=0pt\vbox{\kern #1\box0\kern #2}}

\def\mimic#1#2{\setbox1=\hbox{#1}\setbox2=\hbox{#2}\ht2=\ht1\box2}

\def\stack#1#2#3{\vbox{\offinterlineskip
  \setbox2=\hbox{#2}
  \setbox3=\hbox{#3}
  \dimen0=\ifdim\wd2>\wd3\wd2\else\wd3\fi
  \hbox to \dimen0{\hss\box2\hss}
  \kern #1
  \hbox to \dimen0{\hss\box3\hss}}}

\def\hexp#1{%
  \setbox0=\hbox{${}^{#1}$}%
  \hbox to .5\wd0{\box0\hss}}

\long\def\commentout#1{}

\renewcommand{\eprint}[1]{#1}

\makeatletter

\def\bib@div@mark#1{%
    \@mkboth{{#1}}{{#1}}%
	}
\def\print@backrefs#1{%
    \space\SentenceSpace$\leftarrow$\csname br@#1\endcsname
}
\catcode`\'=11 
\renewcommand{\PrintAuthors}[1]{%
    \ifx\previous@primary\current@primary
        \sameauthors\@empty
    \else
        \def\current@bibfield{\bib'author}%
		        \PrintNames{}{}{\scshape #1}%
   \fi
}
\catcode`\'=12
\def\MRhref#1#2{%
    \begingroup
        \parse@MR#1 ()\@empty\@nil%
        \href{\MR@url}{\texttt{\@tempd\vphantom{()}}}%
        \ifx\@tempe\@empty
        \else
            \ \href{\MR@url}{\texttt{(\@tempe)}}%
        \fi
    \endgroup
}%
\def\MR#1{%
    \relax\ifhmode\unskip\spacefactor3000 \space\fi
    \begingroup
        \strip@MRprefix#1\@nil
        \edef\@tempa{\@nx\MRhref{MR\@tempa}{\@tempa}}%
    \@xp\endgroup
    \@tempa
}
\makeatother

\newcommand\capdraw[5]{
\begin{figure}[tb]
  \setdrawbox{#3}{\artfolder #4}{#5}
  \centerline{\ifgridlines\boxgrid\drawing\fi\box\drawing}
  \caption{#2}
  \label{#1}
\end{figure}%
}

\newarrow{Injectto}{>}{-}{-}{-}{>}
\newarrow{Identifyto}{=}{=}{=}{=}{=}
\newarrow{Rationalto}{}{dash}{}{dash}{>}
\newarrow{IMapsto}{|}{-}{-}{-}{>}

\newdiagramgrid{ilgrid}{0.385, 0.385}{}
\def\ilto{\pdfsyncstop\mathop{\hbox{\kern -5pt
\begin{diagram}[midshaft,inline,grid=ilgrid]
&\rTo&\end{diagram}\kern -5pt}}\pdfsyncstart}
\def\ilinto{\pdfsyncstop\mathop{\hbox{\kern -5pt
\begin{diagram}[midshaft,inline,grid=ilgrid]
&\rInjectto&\end{diagram}\kern -5pt}}\pdfsyncstart}
\def\ilmapsto{\pdfsyncstop\mathop{\hbox{\kern -5pt
\begin{diagram}[midshaft,inline,grid=ilgrid]
&\rIMapsto&\end{diagram}\kern -5pt}}\pdfsyncstart}
\newdiagramgrid{ildashgrid}{0.55, 0.55}{}
\def\ildashto{\pdfsyncstop\mathop{\hbox{\kern -5pt
\begin{diagram}[midshaft,inline,grid=ildashgrid]
&\rDashto&\end{diagram}\kern -5pt}}\pdfsyncstart}

\renewcommand{\to}{\pdfsyncstop\ilto\pdfsyncstart}

\newdiagramgrid{esgrid}{0.585, 0.585}{}
\def\esarrowplain{\pdfsyncstop\mathop{\hbox{\kern -3pt
\begin{diagram}[midshaft,inline,grid=esgrid]
&\rTo&\end{diagram}\kern -3pt}}\pdfsyncstart}
\def\esarrowmap#1{\pdfsyncstop\mathop{\hbox{\kern -3pt
\begin{diagram}[midshaft,inline,grid=esgrid]
&\rTo&\end{diagram}\kern -3pt}}^{\kern -3pt #1}\pdfsyncstart}
\def\esarrowinject{\pdfsyncstop\mathop{\hbox{\kern -3pt
\begin{diagram}[midshaft,inline,grid=esgrid]
&\rInjectto&\end{diagram}\kern -3pt}}\pdfsyncstart}
\def\esarrowinjectmap#1{\pdfsyncstop\mathop{\hbox{\kern -3pt
\begin{diagram}[midshaft,inline,grid=esgrid]
&\rInjectto&\end{diagram}\kern -3pt}}^{#1}\pdfsyncstart}
\newdiagramgrid{esgridshort}{0.47, 0.47}{}
\def\shortesarrowplain{\pdfsyncstop\mathop{\hbox{\kern -3pt
\begin{diagram}[midshaft,inline,grid=esgridshort]
&\rTo&\end{diagram}\kern -3pt}}\pdfsyncstart}
\def\shortesarrowmap#1{\pdfsyncstop\mathop{\hbox{\kern -3pt
\begin{diagram}[midshaft,inline,grid=esgridshort]
&\rTo&\end{diagram}\kern -3pt}}^{\kern -3pt #1}\pdfsyncstart}
\def\shortesarrowinject{\pdfsyncstop\mathop{\hbox{\kern -3pt
\begin{diagram}[midshaft,inline,grid=esgridshort]
&\rInjectto&\end{diagram}\kern -3pt}}\pdfsyncstart}
\def\shortesarrowinjectmap#1{\pdfsyncstop\mathop{\hbox{\kern -3pt
\begin{diagram}[midshaft,inline,grid=esgridshort]
&\rInjectto&\end{diagram}\kern -3pt}}^{#1}\pdfsyncstart}

\newcommand{\setdiagram}[3]{
\ifx#1\endA\endB
\setbox0=\hbox{}
\else
\setbox0=\hbox{(\theequation)}
\fi
\pdfsyncstop
\ifx#2\endA\endB
\begin{diagram}[LaTeXeqno,midshaft]
\label{#1}
#3
\end{diagram}
\else
\begin{diagram}[LaTeXeqno,midshaft,#2]
\label{#1}
#3
\end{diagram}
\fi
\pdfsyncstart
}

\newcommand{\inlinechoose}[2]%
{\smash{\binom{#1}{#2}}}

\newcommand{\inlinefrac}[2]%
{\smash{\bigl({\frac{#1}{#2}} \bigr)}}

\newcommand{\inlinefracnoparen}[2]%
{\smash{ {\frac{#1}{#2}} }}

\newcommand{\CC}{\ensuremath{\mathbb{C}}}

\newcommand{\PP}{\ensuremath{\mathbb{P}}}
\newcommand{\PNminusone}{\ensuremath{\mathbb{P}^{N-1}}}
\renewcommand{\AA}{\ensuremath{\mathbb{A}}} 

\newcommand{\quotient}{{/\!\!/}}
\newcommand{\spec}{\mathrm{Spec}}

\newcommand{\pic}{\mathrm{Pic}}

\newcommand{\codim}{\mathrm{codim}}
\newcommand{\pgl}{\mathrm{PGL}}

\newcommand{\ssl}{\mathrm{SL}}

\newcommand{\diag}{\mathrm{diag}}

\newcommand{\restricted}{\mathrm{res}}

\newcommand{\shavedast}{\ast \kern -1.25pt}

\renewcommand{\epsilon}{\varepsilon}
\newcommand{\pst}{Potential Stability Theorem\xspace}
\newcommand{\plucker}{Pl\"ucker\xspace}
\newcommand{\resp}[1]{[resp: {#1}]}

\newcommand{\thst}[2]{\ensuremath{{#1}^{\mathrm{#2}}}}

 \newcommand{\restrictedto}[2]{%
 \ensuremath{{#1\hskip 1pt{\vrule height 7.2pt depth 3.6pt}\hskip0.75pt}%
 \raisebox{-1.8pt}[0pt][0pt]{\ensuremath{\null_{#2}}}}}

\newcommand{\h}{\ensuremath{\mathbf{H}}}
\newcommand{\hhat}{\ensuremath{\hat{\h}}}
\newcommand{\m}{\ensuremath{M}}
\newcommand{\g}{\ensuremath{G}}
\newcommand{\p}{\ensuremath{P}}
\newcommand{\sss}{\ensuremath{S}}

\newcommand{\M}{\ensuremath{\mathcal{M}}}

\renewcommand{\o}{\ensuremath{\mathcal{O}}} 
\renewcommand{\i}{\ensuremath{\mathcal{I}}} 

\renewcommand{\phi}{\varphi}

\newcommand{\mi}[1]{\ensuremath{\m_{#1}}}
\newcommand{\mg}{\ensuremath{\mi{g}}}

\newcommand{\mij}[2]{\ensuremath{\m_{{#1},{#2}}}}
\newcommand{\mgn}{\ensuremath{\mij{g}{n}}}

\newcommand{\mbar}{\ensuremath{\overline{\m}}}
\newcommand{\mibar}[1]{\ensuremath{\mbar_{#1}}}
\newcommand{\mijbar}[2]{\ensuremath{\mbar_{{#1},{#2}}}}
\newcommand{\mgbar}{\ensuremath{\mibar{g}}}
\newcommand{\mgnbar}{\ensuremath{\mijbar{g}{n}}}

\newcommand{\Mbar}{\ensuremath{\overline{\M}}}
\newcommand{\Mibar}[1]{\ensuremath{\Mbar_{#1}}}
\newcommand{\Mijbar}[2]{\ensuremath{\Mbar_{{#1},{#2}}}}
\newcommand{\Mgbar}{\ensuremath{\Mibar{g}}}

\newcommand{\dhat}{{\hat d}}

\newcommand{\mhat}{{\hat m}}
\newcommand{\Mhat}{{\hat M}}
\newcommand{\uhat}{{\hat u}}
\newcommand{\Uhat}{{\hat U}}
\newcommand{\Vhat}{{\hat V}}
\newcommand{\What}{{\hat W}}
 
\newcommand{\mgnbardhat}{\mgnbar(\PP^r,\dhat)}

\renewcommand{\k}{\ensuremath{\mathcal{K}}}

\newcommand{\defaultarray}{\renewcommand{\arraystretch}{1.7}}

\newcommand{\stretcharray}[3]{\renewcommand{\arraystretch}{#1}%
{{#3}}\defaultarray}

\renewcommand{\hbar}{\ensuremath{\overline{H}}}

\newcommand{\proj}{\mathrm{Proj}}
\newcommand{\mcspan}{\mathrm{span}}

\newcommand{\semi}[1]{\ensuremath{{#1}^{\text{ss}}}}
\newcommand{\stable}[1]{\ensuremath{{#1}^{\text{s}}}}

\newcommand{\sym}{\mathrm{Sym}} 

\newcommand{\ired}[1]{\ensuremath{{#1}_{\mathrm{red}}}}

\newcommand{\normali}[1]{\ensuremath{n_{#1}}} 
\newcommand{\normaly}{\normali{Y}}

\newcommand{\di}[1]{\ensuremath{d_{#1}}}
\newcommand{\dy}{\di{Y}}

\newcommand{\ord}{\mathrm{ord}}

\newcommand{\omegan}[2]{\omega_{#1}^{\otimes #2}}
\newcommand{\omeganu}[1]{\omegan{#1}{\nu}}
\newcommand{\Homeganu}[1]{H^0(#1,\omeganu{#1})}
\newcommand{\homeganu}[1]{h^0(#1,\omeganu{#1})}

\newcommand{\ixon}[4]{\ensuremath{{#1}^{#2}\bigl({#3},\o_{#3}({#4})\bigr)}}
\newcommand{\hixon}[3]{\ixon{h}{#1}{#2}{#3}}
\newcommand{\hzeroxon}[2]{\hixon{0}{#1}{#2}}

\newcommand{\hzeroxom}[1]{\hzeroxon{#1}{m}}

\newcommand{\Hixon}[3]{\ixon{H}{#1}{#2}{#3}}
\newcommand{\Hzeroxon}[2]{\Hixon{0}{#1}{#2}}

\newcommand{\Hzeroxom}[1]{\Hzeroxon{#1}{m}}

\newcommand{\Hoxom}{\Hzeroxom{X}}

\newcommand{\Hocom}{\Hzeroxom{X}}

\newcommand{\invv}[2]{{#1}^{\kern #2{\null^{-1}}}}

\newcommand{\PV}{\PP(V)}
\newcommand{\prodin}{\prod_{i=1}^{n}}

\newcommand{\GIT}{GIT\xspace}

\newcommand{\hvp}{\ensuremath{\h_{\PV,P}}}

\newcommand{\ddiv}{\mathrm{Div}}
\newcommand{\chow}{\mathrm{Chow}}

\newcommand{\lbpow}[2]{#1^{\otimes #2}}
\newcommand{\lbm}[1]{#1^{\otimes m}}

\newcommand{\HO}[2]{H^0\bigl(#1,#2\bigr)}
\newcommand{\hO}[2]{h^0\bigl(#1,#2\bigr)}

\newcommand{\mpsgbar}{\overline{M}_g^{\text{ps}}}

\newcommand{\mhgbar}{\overline{M}_g^{\text{h}}}

\newcommand{\ct}{\widetilde{c}}

\newcommand{\cmin}{\underline{c}}
\newcommand{\cmax}{\overline{c}}

\newcommand{\desiredspread}{1.07}
\makeatletter
\newcommand{\normalspread}{\linespread{\desiredspread}\thm@bodyfont{\itshape}\selectfont}
\makeatother
\newcommand{\beginS}[1]{\begin{#1}\normalspread}
\newcommand{\beginSS}[2]{\begin{#1}{#2}\normalspread}
\renewcommand{\beginS}[1]{\begin{#1}}
\renewcommand{\beginSS}[2]{\begin{#1}{#2}}

\newdiagramgrid{flipgrid}{0.5,1.,0.5,1.,0.5}{0.5,1.,0.5}

\newcommand{\res}[2]{\tau_{#1{\rightarrow}#2}}

\newcommand{\flushsubsection}[1]{\subsection{\kern-\normalparindent #1}}
\newcommand{\flushsubsectionstar}[1]{\subsection*{\kern-\normalparindent #1}}

\begin{document}

\ifstandalone

\title{GIT Constructions of Moduli Spaces of Stable Curves and Maps}

\author{Ian Morrison}
\address{Department of Mathematics\\ Fordham University\\ Bronx, NY 10458}
\email{morrison@fordham.edu}


\begin{abstract}
Gieseker's plan for using GIT to construct the moduli spaces of stable curves, now over 30 years old, has recently been extended to moduli spaces of pointed stable curves and of stable maps by Swinarski and Baldwin. The extensions turn out to be surprisingly delicate and both require the development of novel techniques for checking stability of Hilbert points. Simultaneously, interest in the area has been spurred by the log minimal model program of Hassett and his coworkers Hyeon and Lee in which these models are produced by suitably modified GIT constructions. Here I first give an introduction to the area by sketching Gieseker's strategy. Then I review a number of  variants---those involving unpointed curves that arise in Hassett's program emphasizing Schubert's moduli space of pseudostable curves, that of Swinarski for weighted pointed stable curves, and that of Baldwin and Swinarski for pointed stable maps---focusing on the steps at which new ideas are needed. Finally, I list open problems in the area, particularly some arising in the log minimal model program that seem inaccessible to current techniques.   
\end{abstract}  

\maketitle

\else 
\title{Geometry of Riemann surfaces and their moduli spaces}
\tableofcontents
\mainmatter

\chapter{GIT Constructions of Moduli Spaces of Stable Curves and Maps}

\author{Ian Morrison}
\address{Department of Mathematics\\ Fordham University\\ Bronx, NY 10458}
\email{morrison@fordham.edu}


\begin{abstract}
Gieseker's plan for using GIT to construct the moduli spaces of stable curves, now over 30 years old, has recently been extended to moduli spaces of pointed stable curves and of stable maps by Swinarski and Baldwin. The extensions turn out to be surprisingly delicate and both require the development of novel techniques for checking stability of Hilbert points. Simultaneously, interest in the area has been spurred by the log minimal model program of Hassett and his coworkers Hyeon and Lee in which these models are produced by suitably modified GIT constructions. Here I first give an introduction to the area by sketching Gieseker's strategy. Then I review a number of  variants---those involving unpointed curves that arise in Hassett's program emphasizing Schubert's moduli space of pseudostable curves, that of Swinarski for weighted pointed stable curves, and that of Baldwin and Swinarski for pointed stable maps---focusing on the steps at which new ideas are needed. Finally, I list open problems in the area, particularly some arising in the log minimal model program that seem inaccessible to current techniques.   
\end{abstract}  
\fi

\linespread{\desiredspread}\normalfont\selectfont

\section{Introduction}
\label{intro}

Gieseker first used GIT to construct the moduli space $\mgbar$ of stable curves over 30 years ago. I learned his ideas in writing up Mumford's Fields Medalist lectures~\cite{MumfordEnseignement}, in which $\mgbar$ is realized as a quotient of a suitable Chow variety. Gieseker himself later wrote up versions based on lectures at the Tata Institute~\cite{GiesekerTata} and, later, at CIME~\cite{GiesekerCIME}. In both of these, Hilbert schemes serve as the parameter space and this variant has now become standard. 

The strategy of Gieseker's construction has recently been extended to give GIT constructions of other moduli spaces of stable curves and maps. Even where other constructions of these spaces were known, these GIT constructions are of interest because they come equipped with natural ample classes that can be readily expressed in terms of standard line bundles and divisors. In other cases, these constructions yield new birational models that turn out to arise naturally in running the log minimal model program for these spaces. The aim of this article is to review this work and point out some interesting open problems in the area. 

This introduction gives an informal overview of the main stages of the constructions.  In it, I assume familiarity with the basic GIT setup of Hilbert stability problems, but the unfamiliar reader will find the definitions and results involved in Section~\ref{numerical} and further details about the steps below can be found in Section~\ref{gieseker}. Here it suffices to identify four main steps.

\beginS{enumerate}
	\item Show that Hilbert points of smooth objects embedded by sufficiently ample linear series are GIT asymptotically stable. 
	\item Prove a Potential Stability theorem: that is, show that the Hilbert point of any object embedded by a sufficiently ample linear series can be GIT asymptotically semi-stable only if the object is abstractly or moduli stable in a suitable sense and if, in the reducible case, the components are embedded in a sufficiently balanced way. 
	\item Show that the locus of Hilbert points of $\nu$-canonically embedded objects is, if $\nu$ is large enough, locally closed and smooth or nearly so. 
	\item Show that any $\nu$-canonically embedded moduli stable object not ruled out by the Potential Stability theorem must be Hilbert semi-stable.
\end{enumerate}

This plan has recently been carried out to give GIT constructions of moduli spaces of pointed stable curves and of stable maps by Swinarski~\cite{SwinarskiThesis} and Baldwin~\cite{BaldwinSwinarski}. That the former, at least, of these had not been undertaken long ago is surprising, although Pandharipande in the eprint~\cite{PandharipandeConfiguration} did give a GIT construction of a different compactification of $\mgn$ a decade ago. You'd expect that adding marked points would require only minor modifications of the arguments. In fact, even this extension, recently completed in the thesis of Dave Swinarski~\cite{SwinarskiThesis}, turns out to be rather tricky. Likewise, you'd expect that the step from (pointed) curves to (pointed) maps would pose more serious challenges. Again things turn out unexpectedly. Swinarski~\cite{SwinarskiUnpointed} had earlier constructed moduli spaces of maps from unpointed curves using arguments very close to those of Gieseker. But once more, adding marked points makes checking Hilbert stability much more delicate. A detailed discussion of these difficulties must wait until I have reviewed the construction of $\mgbar$ in the next section and introduced the notions involved. Here, for those with some familiarity with these constructions, I sketch two main ones. 

First, the GIT problems that arise involve a choice of linearization, unlike the unpointed case where the linearization is  canonical up to scaling. The extra parameters on which the linearization now depends must be selected carefully to obtain the desired quotient. In particular, there are some additional technical difficulties in the second step above.  An interesting open problem is to better understand this VGIT problem. Do other linearizations lead to quotients that are moduli spaces for variant moduli problems? More generally, can we describe the VGIT chamber structure of these problems and understand the wall crossing modifications in terms of natural classes on these moduli spaces? How do these variations fit into the log minimal model program for these spaces? For more on these questions, see Section~\ref{open}. 

More seriously, Gieseker's techniques for showing that smooth objects with sufficiently ample polarization are asymptotically  Hilbert stable fail when there are marked points and new ideas are needed. For a scheme $X$ in $\PP^{N-1}$, this involves verifying a numerical criterion (Proposition~\ref{NumericalCriterionForHilbertPoints} is a model) for an arbitrary non-trivial $1$-parameter subgroup $\rho$ of $\pgl(N)$ and that this in turn amounts to showing that $\Hoxom$ has a basis of negative $\rho$ weight for sufficiently large $m$. For a curve $C$ with marked points $p_i$, a contribution from a section not vanishing at each $p_i$ is added to the $\Hocom$ term and the estimates for this latter coming from Gieseker's Criterion~\ref{GiesekersCriterion} are not sharp enough to incorporate the former.

Section~\ref{swinarski} covers Swinarski's construction of moduli of weighted pointed curves. His approach is a refinement of Gieseker's. The idea in both cases is first to exhibit filtrations of $\Hocom$ by subspaces whose weights are bounded and whose dimensions can be estimated by Riemann-Roch and then to verify the numerical criterion by combinatorial arguments using this data. Gieseker works only with \emph{monomial} subspaces, by which I mean subspaces. spanned by monomials in a basis of $\PP^{N-1}$ compatible with the action of $\rho$. Swinarski uses subspaces that are spans of several such monomial subspaces and that I'll call \emph{polynomial}. One consequence of the use of these more complicated subspaces is that the calculations needed to show that the corresponding bases have negative weight become much more delicate because it is not known how to reduce this to a linear programming problem.

Baldwin deals with pointed stable maps by a radically different strategy that is outlined in Section~\ref{baldwin}. She very cleverly relates the numerical criterion for the Hilbert point of a map $f$ in which the underlying curve is $C$ is smooth with respect to a $1$-ps $\rho$ to that for a map $f'$ in which the underlying curve $C'$ is $C$ with an elliptic tail added at the $\thst{n}{th}$ marked point with respect to a related $1$-ps $\rho'$. This idea seems shocking at first, because the projective spaces that are the targets of $f$ and $f'$ are different and that $\rho$ and $\rho'$ lie in different groups. The choices of $f'$ and  $\rho'$ are \emph{not} canonical and the ambiguity can only be resolved pointwise in both cases. Nonetheless, Baldwin is able to show that if $f$ were unstable with respect to  $\rho$ then any $f'$ would necessarily be unstable with respect to  $\rho'$. This makes possible a diagonal induction to the case of unpointed  maps of genus $g+n$ which is treated by Swinarski in~\cite{SwinarskiUnpointed}.

The cases of weighted pointed curves and of pointed maps overlap. The moduli space $\mgnbar$ is the former with all weights equal to $1$ and the latter with target space a point. However, the two approaches are both of interest even in this common case because they handle different sets of GIT problems. Swinarski is able to deal with lower values of the canonical multiple $\nu$ and, for those $\nu$ to which both apply, the sets of linearizations that can be handled do not nest in either direction. 

Simultaneously, interest in the area has been spurred by the log minimal model program of Hassett and his coworkers Hyeon and Lee in which these models are produced by suitably modified GIT constructions. The key construction in the paper~\cite{HassettHyeonLogCanonical} of Hassett and Hyeon is a result of Schubert~\cite{Schubert}. Gieseker's construction of $\mgbar$ requires taking the canonical multiple $\nu \ge 5$. Schubert worked out what happens if we take $\nu=3$ and shows that the resulting Chow quotient is a moduli space $\mpsgbar$ for what he dubs \emph{pseudostable} curves on which ordinary cusps are allowed but elliptic tails are not. His construction is able to hew closely to Gieseker's except, of course, at the points in steps 2. and 4. at which curves with cusps and with elliptic tails are handled and is easily modified to see that the Hilbert quotient is again $\mpsgbar$. Although Schubert does not treat $\nu=4$, his proof was widely assumed to apply also to this case. Recently, Hyeon and the author~\cite{HyeonMorrison} were led to examine this assumption and discovered that some additional refinements are needed but that both the Chow and Hilbert quotients are again $\mpsgbar$. These constructions are discussed in Section~\ref{schubert} which also introduces some new ideas that arise in more recent work of Hassett and Hyeon~\cite{HassettHyeonFlip} on the case $\nu=2$ and points to Hassett's study~\cite{HassettGenusTwo} of the genus $2$ case that launched work in this area.

Finally, section~\ref{open} discusses open problems in this area. These fall into two main groups. Possible streamlinings of some of the recent constructions and the VGIT problems mentioned above form one group. The second involves GIT problems that arise out of the log minimal model program. These ask for descriptions, for pointed curves, of the quotients that result from using a $\nu$ smaller than that required to produce the moduli space, both intrinsically as moduli spaces for a variant moduli problem \`a la Schubert and as log models. To introduce these, I review one small genus example of Hyeon and Lee~\cite{HyeonLeeGenusThree} and point to related work of Smyth~\cite{SmythThesis}.

However, it also appears that, as we approach the canonical models of these spaces, more delicate questions arise. The specification of a Hilbert stability problem involves not only the choice of $\nu$ but also that of the ``sufficiently large'' degree $m$. I review calculations of Hassett based on the results of~\cite{GibneyKeelMorrison} that predict that log minimal models with scaling for $\mgbar$ arise as the quotients that result for small values of $\nu$ and \emph{fixed} values of $m$. Answering these questions will require completely new ideas since all existing techniques for checking Hilbert stability prove this only asymptotically and hence require $ m \gg 0$. 

\flushsubsectionstar{Acknowledgements} Since this paper is essentially expository, I owe a great deal to the authors whose ideas I have tried to explain. My own contribution has been limited to trying to clarify and simplify arguments where I could. Readers may judge with what success when they return to the primary sources to fill in the many steps I had to omit here. I have also tried to balance the mutually exclusive aims of having the notation be internally consistent and be consistent with these sources. My rule was to give priority to the former while trying to stay close to the latter and to supply dictionaries when the two strayed too far apart.

Much of the writing of this survey was completed in the spring of 2008 while I was visiting the University of Sydney with support from a Fordham University Faculty Fellowship. In addition to thanking both these institutions, I would like to express my gratitude for the hospitality shown to me while I was in Australia by Gus Lehrer, Amnon Neeman and Paul Norbury.

I have also benefitted from discussions of and correspondence about the ideas discussed here with Elizabeth Baldwin, Dave Bayer, Joe Harris, Brendan Hassett, Julius Ross and Michael Thaddeus. Finally, a special thank you to the Davids, Gieseker and Mumford who taught me the subject and Hyeon and Swinarski who explained not only their recent work but the questions that arise from it. 

 \section{Stability of Hilbert points} \label{numerical}

\subsection{Setup and Linearization}
\label{setupandlinearization}

This goal of this section is to understand the numerical criterion for the $\pgl(V)$-action on the Hilbert scheme $ \h=\hvp$ of subschemes of $\PV \cong \PNminusone$ with Hilbert polynomial $\p(m)$. The ``$-1$'' is inserted above because it is then $N$ and not, as in most sources, $N+1$ that will be ubiquitous in later formulae. The same goal prompts the unusual indexing, starting at $1$, of homogeneous coordinates that will appear shortly. With these, and other notational changes made to conform with the notation of the constructions I will be summarizing, the treatment here follows closely that in \cite{Moduli}*{Section 4.B} to which the reader is referred for further details. See~\cite{Moduli}*{Section 4.A} for a quick review of more basic notions in \GIT\ and~\cite{GIT} for a thorough one. 

We need to first recall the procedure for linearizing this action. Recall that this means lifting the action of $\pgl(V)$ to one on an ample line bundle $L$ on $\h$ and, in turn on the sections of $L$. It is convenient to pass first to the finite cover $\ssl(V)$---and harmless since the scalar matrices corresponding to $\thst{N}{th}$ roots of unity act trivially on $\PV$. The action of $\ssl(V)$ on $\PV$ then lifts to its natural action on $V$ and hence yields an linearization on the line bundle $\o_{\PV}(1)$. It's convenient (and hence standard) to express the numerical criterion in terms of the $\ssl(V)$ action.

Further, fixing a sufficiently large degree $m$, the action of $\ssl(V)$ on $V$ induces, in turn, actions on $S_m := \sym^{m}(V)^{{\vee}}$ and on $W_m=\bigwedge^{\p(m)}\big(S_m\big)$. In the same way, that the action of $\ssl(V)$ on $V$ gives a linearization of $\o_{\PV}(1)$, its action on $W_m$ gives a linearization on $\o_{\PP(W_m)}(1)$. The Hilbert scheme $\h$ has a natural Pl\"ucker embedding in $\PP(W_m)$ as a subscheme of the Grassmannian \g\ of $P(m)$-dimensional quotients $Q$ of $S_m$. Under this embedding, $\o_{\PP(W_m)}(1)$ restricts to the tautological very ample invertible sheaf $\Lambda_m$ on $\h$ that thus also acquires an $\ssl(V)$-linearization depending only on the choice of the degree $m$. An equivalent description of $\Lambda_m$ is as $ \det(\pi_*(\o_X(m))$ where $\pi:X \to \h$ is the universal family and $\o_X(1)$ is the tautological polarization: see~\cite{HassettHyeonFlip}*{Proposition 3.10} for details. 

\subsection{The numerical criterion for Hilbert points}
\label{numericalcriterion}

On, then, to the numerical criterion (see~\cite{Moduli}*{Section 4.A}). We fix a one-parameter subgroup $ \rho:\CC^{*}\to\ssl(V)$ and homogeneous coordinates $B=B_{\rho}=\{x_{1},\ldots,x_i, \ldots, x_{N}\} $ that we view as a basis of $V^{\vee}$ with respect to which 
\beginS{equation}\label{setuponeps}
 \rho(t)=\diag\left(t^{w_{1}},\ldots,t^{w_{i}}, \ldots,t^{w_{N}}\right) 
\end{equation}
with $\sum_{i}^{N}w_{i}=0$. The data of $\rho$ is thus equivalent to the data of $B$ considered as a \emph{weighted basis} (i.e., along with a set of integral weights $w_{i}$ summing to 0) and we'll henceforth refer to $B$ and $\rho$ interchangeably.

The weighted basis $B_m$ consisting of degree~$m$ monomials $y = \prodin x_i^{m_i}$ in the $x_{i}$'s with weights $w(y)=\sum_{i}^{N}w_{i}m_{i}$ diagonalizes the action of $\rho$ on $S_m$. Likewise, the \plucker\ basis consisting of all unordered $\p(m)$-element subsets $ z=\{y_{j_{1}},\ldots,y_{j_{\p(m)}}\}$ of $B_{m}$ with weights $w(z) := \sum_{k=1}^{\p(m)} w(y_{j_k})$ diagonalizes the action of $\rho$ on $W_m$.

The key observation is that the \plucker\ coordinate $z$ is nonzero at the point $[Q]$ of the Grassmannian \g\ corresponding to a quotient $Q$ if and only if the images in $Q$ of the $\p(m)$ monomials in $z$ form a basis of $Q$. Since the Hilbert point $[X]$ of a subscheme of $X$ of $\PV$ with Hilbert polynomial $\p(m)$ corresponds to the quotient
$$
S_m = \sym^{m}(V)^{{\vee}} \rTo^{~~\restricted_{X}~~} \Hzeroxom{X}\, ,
$$
$Z$ is nonzero at $[X]$ if and only if the restrictions $ \restricted_{X}(y_{j_{k}})$ of the monomials in $z$ are a basis of $\Hzeroxom{X}$. We will call such a set of monomials a \emph{$B$-monomial basis} of $\Hzeroxom{X}$.

Before using this observation to interpret the numerical criterion, it is convenient to push the change of point of view from one-parameter subgroups $\rho$ to weighted bases $B$ little further. First, note that, in the language of weighted bases, there is no need to maintain the requirements that the weights $w_{i}$ be integral or sum to 0. Instead, we denote this sum by $w_{B}$. 

The second simplification involves the notion of a rational \emph{weighted filtration} $F$ of $V$. This is just a collection of subspaces $U_{w}$ of $V$, indexed by the rational numbers, with the property that $ U_{w}\subset U_{w'}$ if and only if $w\geq w'$. Any weighted basis $B$ determines a weighted filtration $F_{B}$ by taking $ U_{r}=\mcspan\{x_{i}|w_{i}\leq w\}$. We say that $B$ is compatible with $F$ if $F_{B}=F$. If so, then we define the weight $w_{F}$ of $F$ to be $w_{B} $: this clearly doesn't depend on which compatible $B$ we choose.

Each $F$ is determined by the subspaces associated to the finite number of $w$ at which there is a jump in the dimension of $U_{w}$. It's convenient to use a notation that implicitly assumes that all these jumps in dimension are of size 1 and to view $F$ as the collection of data: 
{\stretcharray{1.2}{0pt}{
\beginS{equation}
\label{FiltrationNotation}
{\begin{array}{rllllllll}
F=F_{1}:\,V=&V_{1}&\supsetneqq&V_{2}& \supsetneqq&\cdots&\supsetneqq &V_{N}&\supsetneqq\{0\}\\
&w_{1}&\geq&w_{2}&\geq&\cdots&\geq&w_{N}& \end{array}}
\end{equation}%
}%
Thus, $U_{r}=\cup_{w_{i}\leq w} V_{i}$ and an element $x$ in $ V$ has weight $w(x)=w_{i}$ if and only if $x$ lies in $V_{i} $ but not in $V_{i+1}$. Of course, whenever $w_{i}=w_{i+1}$, then $ F$ has a larger jump and $V_{i+1}$ is neither uniquely determined by $F$ nor, indeed, needed to recover the filtration $F$. This harmless ambiguity makes it possible to use the same indexing in discussing one-parameter subgroups, weighted bases and weighted filtrations.

By repeating the arguments above using any basis $B$ compatible with $F$, we see that $F$ determines weighted filtrations $F_{m}$ of each $S_m = \sym^{m}(V)^{{\vee}}$. But anytime we have a weighted filtration on a space $S$ and a surjective homomorphism $\phi:S\to H$, we get a weighted filtration on $H$ by the rule that the weight of an element $h$ of $H$ is the minimum of the weights of its preimages in $S$. Thus, $ F_{m}$ determines by restriction to $X$ a weighted filtration, that we also denote by $ F_{m}$, on \Hzeroxom{X}. We let $w_{F}(m)$ denote the weight of any basis of \Hzeroxom{X} compatible with the filtration $F_{m}$: as the notation suggests, we'll shortly be viewing these weights as giving a function of $m$ depending on $F$. With these preliminaries, we have:

\beginSS{proposition}{[Numerical criterion for Hilbert points]} \label{NumericalCriterionForHilbertPoints} 
The \thst{m}{th} Hilbert point $[X]_{m}$ of a subvariety $X$ of $\PV$ with Hilbert polynomial $\p$ is stable \resp{semistable} with respect to the natural $\ssl(V)$-action if and only if the equivalent conditions below hold:
\beginS{enumerate}
	\item For every weighted basis $B$ of $V$, there is a $B$-monomial basis of \\ $\Hzeroxom{X}$ whose $B $-weights have negative \resp{nonpositive} sum.
	\item For every weighted filtration $F$ of $V$ whose weights $ w_{i}$ have average~$\alpha$, $$ w_{F}(m)< \text{\resp{$\le$}~} m\alpha\p(m)\, . $$
	\end{enumerate}
\end{proposition}

\beginS{proof} The first statement is an immediate translation of the Numerical Criterion~\cite{GIT}*{Theorem~2.1}: if we diagonalize the action of the one-parameter subgroup $\rho$ associated to $B$ on $W$ as above, then the $B$-monomial bases are just the nonzero \plucker\ coordinates of $[X]_{m}$ and their weights are the weights of $[X]_{m}$ with respect to $\rho$. In other words, the Hilbert-Mumford index $\mu_{\rho}([X], \Lambda_m)$ whose sign determines the stability of $[X]$ with respect to $\rho$ and the linearization $\Lambda_m$ equals the least weight $w_B(m)$ of a $B$-monomial basis of $\Hzeroxom{X}$.

To see 2), observe that if $B$ is any basis compatible with the filtration $F$ and we set $w'_{i}=\beta(w_{i}-\alpha)$ where $\beta$ is chosen so that all the weights $w'_{i}$ are integral, then $B$ becomes a weighted basis, and, moreover, every weighted basis $B$ arises in this way from some $F$. The $F$-weight of any degree~$m$ monomial then differs from its $B$-weight by $m\alpha\beta$. Hence the weight of any $B$-monomial basis of \Hzeroxom{X} will differ from $\beta w_{F}(m)$ by $\beta m\alpha \hzeroxom{X}=\beta m\alpha\p(m)$. Therefore, the given inequality is equivalent to the negativity of the $B$-weights of such bases.
\end{proof}

\textsc{Notational remark~} Because all our verifications of stability and instability involve estimating weights of bases, we have stated this (and variant numerical criteria that follow) in such terms. Since, we will always be working with fixed choice of linearization, we have, to simplify, omitted this choice from the notation for such weights.  All these criteria have straightforward translations in terms of the Hilbert-Mumford indices that we henceforth leave to the reader.

We will continue to write $\alpha:=\alpha_{F}$ for the average weight of an element of a basis $B$ of $V$ compatible with $F$. We will also say simply that the variety $X$ is \emph{asymptotically Hilbert stable with respect to $F$} if, for all large $m$, the inequalities of the proposition hold for $F$, and that $X$ is \emph{asymptotically Hilbert stable} if for all large $m$, the \thst{m}{th} Hilbert points of $X$ are stable: i.e., the inequalities of the proposition hold for every nontrivial $F$. All the methods of verifying the stability of an \thst{m}{th} Hilbert point that arise here apply to all sufficiently large $m$, the implicit lower bound depending only on the Hilbert polynomial $\p$ of $X$ so this will not introduce any ambiguity. To see why this is so, we introduce an idea developed in \cite{MumfordEnseignement}: the weights $w_{F}(m)$ are given for large $m$ by a numerical polynomial in $m$ of degree~$(\dim(X)+1)$. For our purposes, all we'll need is the:

\beginSS{lemma}{[Asymptotic numerical criterion]} \label{AsymptoticNumericalCriterion} 
Let $X$ be a subscheme of dimension $r$ and degree~$d$ in \PV. 
	\beginS{enumerate} 
		\item \label{mumfordboundedness} There are constants $C$ and $M$ depending only on the Hilbert polynomial $\p$ of $X$, and, for each $F$, a constant $ e_{F}$ depending on $F$ such that, for all $m\geq M$, 
$$ \left|w_{F}(m)-e_{F}\frac{m^{r+1}}{(r+1)!}\right|<Cm^{r}. $$ 
		\item \label{ealphaInequality} If $e_{F}<\alpha_{F}(r+1)d$, then $X$ is Hilbert stable with respect to $F$; and if $e_{F}>\alpha_{F}(r+1)d$, then $ X$ is Hilbert unstable with respect to~$F$. 
		\item \label{asymptoticuniform} Fix a Hilbert polynomial $\p$ and a subscheme $\sss$ of $\h$. Suppose that there is a $\delta>0$ such that $$ e_{F}<\alpha_{F}(r+1)d-\delta $$ for all weighted filtrations $F$ associated to the Hilbert point of any $X$ in $\sss$. Then there is an $M$, depending only on $\sss$, such that the \thst{m}{th} Hilbert point $[X]_{m}$ of $X$ is stable for all $m\geq M$ and all $X$ in $\sss$. 
	\end{enumerate} 
\end{lemma}

\beginS{proof} For the first assertion, due to Mumford, we'll simply refer to~\cite{MumfordEnseignement}*{Theorem 2.9}. The second then follows by taking leading coefficients in the second form of the numerical criterion and using Riemann-Roch to provide the estimate, for large $m$, $ \p(m)=\hzeroxom{X}=\frac{d}{r!}m^{r}+O(m^{r-1})$. This comparison of leading coefficients shows that $ w_{F}(m)$ will be negative for $m$ greater than some large $M $, but exactly how large this $M$ must be taken depends on the ratio of the constant $C$ in part~1) to the difference $ \alpha_{F}(r+1)d-e_{F}$ in part~2). To get the uniform assertion of part~3), we need both a uniform lower bound (given by $\delta$) for this last difference and the uniform upper bound, provided by Mumford, for $C$. \end{proof}
	
Mumford's argument likewise gives a criterion for Chow stability that we'll need to refer to in Section~\ref{schubert}. Since we won't use Chow points to construct moduli spaces, we'll simply quote it. 
\beginS{corollary}\label{chowstabilitycriterion}
If $e_{F}<\alpha_{F}(r+1)d$, then $X$ is Chow stable with respect to $F$; and if $e_{F}>\alpha_{F}(r+1)d$, then $ X$ is Chow unstable with respect to~$F$. 	
\end{corollary}

\subsection{Gieseker's criterion for curves}
\label{giesekerscriterion}

This subsection reviews a fundamental estimate due to Gieseker for $e_{F}$ that is the main tool for proving Hilbert stability for smooth curves. Although it is not sharp enough to yield stability of smooth curves with marked points, both proofs of this fact that we'll review incorporate many of same ideas. Since this is the only case we'll need, I'll simplify by sticking to curves and I'll omit the combinatorics. 

So fix $C$, a smooth curve embedded in $\PV$ by a linear series with a fixed Hilbert polynomial $\p$ and fix a weighted filtration $F$ as in~\eqref{FiltrationNotation} above. We want to estimate $e_{F}$ in terms of its weights $w_j$ and a new set of invariants, the degrees $d_{j}$ of the subsheaves generated by the sections in the sub-linear series $|V_{j}|$. 

Gieseker first fixes a subsequence 
$$
1=j_{0}>j_{1}>\cdots>j_{h}=N
$$
of $(1,\ldots,N)$. He next introduces two auxiliary positive integers $p$ and $n$ to be fixed later, sets and considers the filtration of $ H^{0}\bigl(C,\o_{C}\bigl(n(p+1)\bigr)\bigr)$ given by the images $U^{n}_{k,l}$ under restriction to $C$ of the subspaces 
$$ W^{n}_{k,l}=\sym^{n}\Bigl(V\cdot\sym^{(p-l)}(V_{j_{k}}) \cdot\sym^{l}\bigl(V_{j_{(k+1)}}\bigr)\Bigr) $$
of $\sym^{n(p+1)}(V)$ where the index $k$ runs from 0 to $ h-1$ and, for each $k$, $l$ runs from $0$ to $p$.

Setting $m=n(p+1)$, give a doubly-indexed filtration of $\Hzeroxon{C}{m}$
{\stretcharray{1.2}{-8pt}{
\begin{equation} {\begin{array}{r@{~}l@{~}l@{~}l@{~}l@{~}l@{~}l@{~}l@{~}l@{~}l@{~}l}%
\label{giesekerfiltration}
 \Hzeroxon{C}{m}=&U^n_{0,0} & \supset & U^n_{0,1} & \supset & \cdots & \supset & U^n_{0,p-1} & \supset & U^n_{0,p} & \\ =&U^n_{1,0} & \supset & U^n_{1,1} & \supset & \cdots & \supset & U^n_{1,p-1} & \supset & U^n_{1,p} & \\ =&\cdots&&&&&&&&\\ =&U^n_{h-1,0} & \supset & U^n_{h-1,1} & \supset & \cdots & \supset & U^n_{h-1,p-1} & \supset & U^n_{h-1,p} & \\ =&U^n_{h,0}\,.&&&&&&&&& \end{array}%
}%
\end{equation}
}
}%
Any element of $U^{n}_{k,i}$ has weight at most  $w_{k,i}=n\bigl(w_{0}+(p-i)w_{j_{k}}+iw_{j_{k+1}}\bigr)$ so any basis of $\Hzeroxon{C}{m}$ compatible with this filtration can have weight at most 
\beginS{equation*}
\beginS{split}
&\Bigl(\sum^{h-1}_{k=0}\sum^{p-1}_{i=0}\left(
 \dim(U^{n}_{k,i})-\dim(U^{n}_{k,i+1})\right)w_{k,i}\Bigr) + \dim(U^{n}_{h,0})w_{h,0}\\
&=\dim(U^{n}_{0,0})w_{0,0}+\Bigl(\sum^{h-1}_{k=0}\sum^{p}_{i=1}\dim(U^{n}_{k,i})\left(w_{k,i}-w_{k,i-1}\right)\Bigr).
\end{split}
\end{equation*}

Gieseker's key claim is that, for any fixed choice of Hilbert polynomial $\p$ and integers $n$ and $ p$, there is an $M$ depending only on these three choices but \emph{not} on the Hilbert point $[C]$ or the weighted filtration $ F$ being considered, such that the dimension formula
\beginS{equation}
\label{GiesekerDimensionFormula} \dim(U^{n}_{k,l})=n\left(d+(p-l)d_{j_{k}}+ld_{j_{k+1}}\right)-g+1 
\end{equation}
holds for every $n\geq M$ and for every $k$ and $i$.

To see \eqref{GiesekerDimensionFormula} pointwise, observe that, if $L_{j}$ is the line bundle on $C$ generated by the sections in $V_{j}$, then we can view $U^{n}_{k,l}$ as a sub-linear series of $H^{0}\bigl(C,(M_{k,l})^{\otimes n}\big)$ where $ M_{k,l}=L\otimes(L_{j_{k}})^{(p-l)}\otimes(L_{j_{k+1}})^{l} $. Since $ |V|$ is tautologically a very ample linear series on $C$ and each $L_j$ is, by definition, generated by the sections in $|V_j|$, the subseries  $U^{1}_{k,l}$ of $H^{0}\bigl(C,M_{k,l}\big)$ is very ample and base point free. Hence, for large $n$, $U^{n}_{k,l} = \sym^n\bigl(U^{1}_{k,l}\bigr)$ will be all of $H^{0}\bigl(C,M_{k,l}\big)$
and the claim follows from Riemann-Roch. The uniform version follows by using standard boundedness arguments to show that $n$ can be chosen to depend only on $\p$, $N$ and $p$. 

From this point on, the argument involves purely formal manipulations that I omit. For details, see~\cite{GiesekerGlobal} or~\cite{Moduli}*{4.B}.  Normalizing so that $w_N = 0$, these lead first to the estimate 
\beginS{equation}
\label{singleFestimate}
\beginS{split}
 e_{F} &\leq 2dw_{0}+\sum^{h-1}_{k=0}\bigl(d_{j_{k}}+d_{j_{k+1}}\bigr) \bigl(w_{j_{k+1}}-w_{j_{k}}\bigr)\\
&= \sum^{h-1}_{k=0}\bigl(e_{j_{k}}+e_{j_{k+1}}\bigr) \bigl(w_{j_{k}}-w_{j_{k+1}}\bigr) 
\end{split}
\end{equation}
where $e_{j} :=d-d_{j}$\ so that $e_{j}$ is the \emph{codegree}, or drop in degree, under projection to $|V_{j}|$.

We now take $\epsilon_F$ to be the minimum of the right hand side of \eqref{singleFestimate} over all subsequences of $\{1, \ldots, N\}$. Since $e_F \le \epsilon_F$, the Asymptotic Numerical Criterion~\ref{AsymptoticNumericalCriterion} immediately gives the first assertion in the following lemma. Because the right hand side of \eqref{singleFestimate} increases if we increase any $e_j$, the inequality in the second part implies stability with respect to any non-trivial $F$. The third assertion then follows by applying the uniform version~\ref{AsymptoticNumericalCriterion}.\ref{asymptoticuniform}.

\beginSS{lemma}{[Gieseker's criterion for curves]}\hfill%
\label{GiesekersCriterion}
\beginS{enumerate} 
	\item A curve $C$ is Hilbert stable with respect to a filtration $F$ with $w_{r}=0$ if $\epsilon_{F}<2d\alpha_{F}$. 
	\item Fix a curve $C$ of degree~$d$ and genus~$g$ in $\PV$ as above, and numbers $\epsilon_{i}$ that are upper bounds for the codegree of \emph{every} subspace $V_{i}$ of codimension $i$ in $V$ and let
$$ \qquad \epsilon_{C}= \max_{\scriptscriptstyle w_{1}\geq \cdots\geq w_{N}=0\atop\scriptscriptstyle\sum_{i=1}^{N}w_{i}=1} \Biggl(\min_{\scriptscriptstyle 1=j_{0}<\cdots<j_{h}=N} \Bigl(\sum^{h-1}_{k=0} \bigl(\epsilon_{j_{k}}+\epsilon_{j_{k+1}}\bigr)\bigl( w_{j_{k}}-w_{j_{k+1}}\bigr)\Bigr)\Biggr)  
$$ 
Then, $C$ is Hilbert stable if 
$ \epsilon_{C}<2\frac{d}{N}\, . $%
	 \item Fix integers $d$, $g$ and $N$ and a subscheme $ \sss$ of the Hilbert scheme of curves of arithmetic genus~$g$ and degree~$d$ in $\PP^{N-1}$. If there is a $\delta>0$ such that $ \epsilon_{C}< 2\frac{d}{N}-\delta $ for every curve $C$ in $\sss$, then there is an $M$ such that the \thst{m}{th} Hilbert point $[C]_{m}$ of $C$ is stable for all $m\geq M$ and all curves $C$ in~$\sss$. \end{enumerate} 
\end{lemma}

\subsection{Stability of smooth curves} 
\label{stabilityofsmoothcurves}

We're now ready to tackle the fundamental:
\beginSS{theorem}{[Stability of smooth curves of high degree]} \label{StabilityOfSmoothCurves} 
Suppose that $C$ is a smooth curve of genus $g\geq 2$ embedded in $\PV$ by a complete linear system $L$ of degree $ d\geq 2g+1$. Then $C$ is asymptotically Hilbert stable. Moreover, an $M$ such that the \thst{m}{th} Hilbert point $[C]_{m}$ is stable for all $m\geq M$ may be chosen uniformly for all such curves~$C$. 
\end{theorem}

We will follow the argument given in \cite{MumfordEnseignement}. The only geometric ingredient is the claim that, for some $\delta >0$, we can take $ \epsilon_{j}=\big(\frac{d}{N-1} - \delta\bigr)(j-1)$ in Gieseker's Criterion. This is most easily seen from the graph in Figure~\ref{RiemannRochLineDiagram} in which the Riemann-Roch line $d=N+g-1$ and the Clifford line $d=2(N-1)$ are graphed in the $(d,N)$-plane. 
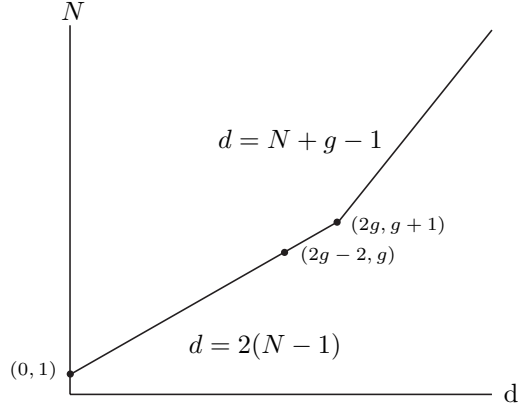
\capdraw{RiemannRochLineDiagram}{Riemann-Roch and Clifford Lines}{.60}{RRClifford.pdf}{
\stext{.045}{.99}{$N$} 
\wtext{.984}{.06}{d} 
\setext{.63}{.16}{$d=2(N-1)$} 
\wtext{.65}{.484}{{\scriptsize$(2g,g+1)$}} 
\wtext{.54}{.408}{{\scriptsize$(2g-2,g)$}} 
\etext{.01}{.12}{{\scriptsize$(0,1)$}} 
\setext{.72}{.67}{$d=N+g-1$}
}%
The corresponding theorems state that the point $(\deg(U),\dim(U))$ corresponding to \emph{any} linear series on $C$ lies in the region below the graph. In particular, this applies to the point $\bigl(d_{j},N-(j-1)\bigr)=\bigl(d-e_{j},N-(j-1)\bigr)$ associated to any linear series $V_{j}$ of codimension $j-1$ in $ H^{0}(C,L)$. On the other hand, the hypothesis of the theorem is that the point $(d,N)$ corresponding to the line bundle $L$ on $C$ lies \emph{on} the Riemann-Roch line. Together, these observations imply that the slope of the line segment from $\bigl(d-e_{j},N-(j-1)\bigr)$ to $(d,N)$ is greater than the slope of the segment joining $(d,N)$ to the ``origin'' $(0,1)$. This is the claim with $\delta=0$, and the claim for small enough positive $\delta$ follows because there are only finitely many choices for the endpoint $\bigl(d-e_{j},N-(j-1)\bigr)$.

Plugging the claim into Gieseker's criterion, we are reduced to checking the following combinatorial claim: 
\beginS{lemma}
\label{combinatorialclaim}
$$\max_{\scriptscriptstyle w_{1}\geq \cdots\geq w_{N}=0\atop\scriptscriptstyle\sum_{i=1}^{N}w_{i}=1} \Biggl(\min_{\scriptscriptstyle 1=j_{0}<\cdots<j_{h}=N} \Bigl(\sum^{h-1}_{k=0} \bigl((j_{k}-1)+(j_{k+1}-1)\bigr)\bigl( w_{j_{k}}-w_{j_{k+1}}\bigr)\Bigr)\Biggr) \le \frac{N-1}{N}$$
\end{lemma}
\beginS{proof}
To carry this out, let's first fix the $w_i$'s. For each subsequence $1=j_{0}>j_{1}>\cdots>j_{h}=N$, consider the ``graph'' obtained by joining the points $ (j_{k}-1,w_{j_{k}})$ and $(j_{k+1}-1, w_{j_{k+1}}) $ by straight line segments as shown in Figure~\ref{GiesekersCriterionSubsequence}. The key observation is that the sum in Gieseker's Criterion corresponding to each subsequence equals twice the area in the first quadrant bounded by the axes and this graph: just integrate with respect to $w$. Taking the minimum of these sums over all  subsequences amounts to computing twice the area under the lower convex envelope $E$ of \emph{all} the points $ (\epsilon_{i},w_{i})$.

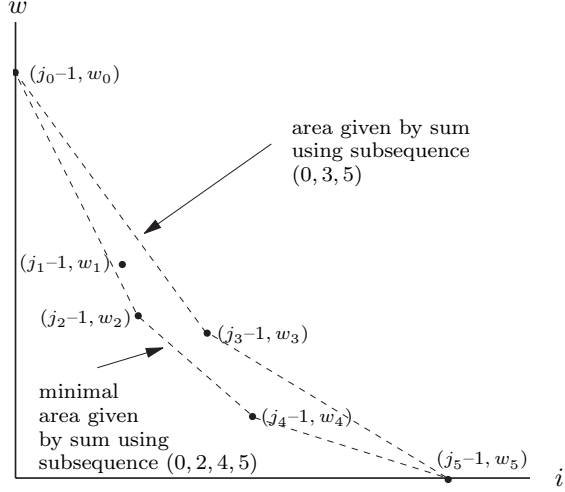
\capdraw{GiesekersCriterionSubsequence}{Area in Gieseker's Criterion}{.60}{GiesekerArea.pdf}{
\swtext{-.01}{1.02}{$w$}
\swtext{.03}{.88}{{\scriptsize$(j_0\mbox{--}1,w_0)$}}
\setext{.19}{.48}{{\scriptsize$(j_1\mbox{--}1,w_1)$}}
\setext{.23}{.365}{{\scriptsize$(j_2\mbox{--}1,w_2)$}}
\swtext{.395}{.33}{{\scriptsize$(j_3\mbox{--}1,w_3)$}}
\swtext{.48}{.155}{{\scriptsize$(j_4\mbox{--}1,w_4)$}}
\swtext{.82}{.065}{{\scriptsize$(j_5\mbox{--}1,w_5)$}}
\wtext{.54}{.78}{{\footnotesize\mimic{subsequence}{area given by sum}}}
\wtext{.54}{.73}{{\footnotesize using subsequence}}
\wtext{.54}{.68}{{\footnotesize\mimic{subsequence}{$(0,3,5)$}}}
\wtext{.05}{.22}{{\footnotesize\mimic{subsequence}{minimal}}}
\wtext{.05}{.17}{{\footnotesize area given}}
\wtext{.05}{.12}{{\footnotesize by sum using}}
\wtext{.05}{.07}{{\footnotesize\mimic{subsequence}{subsequence $(0,2,4,5)$}}}
\swtext{1.05}{.02}{$i$}
}%

Now allow the $w_{i}$'s to vary. If any of the points $(i-1,w_{i})$ does \emph{not} lie on $E$, then moving it down onto $E$ will leave the minimum in Gieseker's Criterion unchanged while reducing the sum of the $w_{i}$'s. Dually, this means that the maximum over sets of weights summing to $1$ in Gieseker's Criterion must occur when the weights are chosen so that \emph{all} the points $ (i-1,w_{i})$ lie on $E$. For such weights, the sum associated to the full sequence---that is, $j_{i}=i$ for all $ i$ from $1$ to $N$---realizes the minimum over all subsequences. If we now compute the area under $E$ by ``integrating with respect to $i$'', we get \beginS{equation*}
\beginS{split}
\quad& \sum_{i=1}^{N-1} \frac{1}{2}(w_i+w_{i+1}) \\
= &\sum_{i=1}^N w_i -\frac{1}{2}(w_1+w_{N})\\
\ge &\sum_{i=1}^N w_i -\frac{1}{N}(\sum_{i=1}^N w_i) \text{\quad by convexity of the weights}\\
= &\frac{N-1}{N} \text{\quad since the weights sum to 1.} \qedhere
\end{split}
\end{equation*}
\end{proof}

\section{Gieseker's Construction of $\mgbar$} 
\label{gieseker}

\subsection{Overview} The goal of this section is to outline the main ideas in Gieseker's \GIT\ construction of \mgbar\ when $ g\geq 2 $. I'll begin with a quick precis for the benefit of those who want to get quickly to newer constructions, then flesh this out for those who are seeing these constructions for the first time. I have given very few proofs in this section since most of the details can be found in~\cite{Moduli}*{Sections 4.B and 4.C}.

The natural approach is to show that suitable pluricanonical models of Deligne-Mumford stable curves have stable Hilbert points and apply \GIT. (It's convenient to use the term pluricanonical even when the curve is singular, understanding that $\omega_C$ is intended where $K_C$ is named.) The first step in such an approach is to define suitable pluricanonical loci and show that they are locally closed in the relevant Hilbert scheme and smoothe (or, in the sequel, nearly so).   

For smooth curves, Hilbert stability of $\nu$-canonical models for $\nu \ge 2$ is immediate from Theorem~\ref{StabilityOfSmoothCurves}. However, no direct proof that Hilbert points of singular Deligne-Mumford stable curves verify the numerical criterion is known. In particular, as shown in~\cite{Moduli}*{Exercise 4.32}, Gieseker's criterion may fail for such points. 

Instead, an indirect approach is used to verify the stability of certain Hilbert points of singular stable curves. The first step is to prove a \pst\ for curves $C$ embedded by an invertible sheaf $L$ degree sufficiently large relative to the arithmetic genus. Such a theorem shows that if such a curve isn't Deligne-Mumford semistable then it has a nonsemistable Hilbert point, and if it's Deligne-Mumford semistable and reducible, then the degree of $L$ on any subcurve $D$ must be approximately proportional to the genus of $D$. For $\mgbar$, the relevant result is Theorem~\ref{PotentialStabilityTheorem}.

The second step of the indirect approach involves considering a smoothing over a discrete valuation ring of a pluricanonically embedded stable curve $ C $. The pluricanonical Hilbert points of the smooth fibers in such a family are stable so, by a semistable replacement argument, we can, after a base change, if necessary, assume that the special fiber is Hilbert semistable. The \pst\ is then used to deduce that this limit can only be the Hilbert point of the pluricanonical model of $ C $. In the construction of $\mgbar$ this step is straightforward. It was extended by Caporaso~\cite{CaporasoCUP}, at the cost of considerably greater technical complications, to prove a converse (requiring somewhat larger $d$) to the \pst\ that she then applied to construct modular compactifications of the universal Picard varieties of degree $d$ line bundles over $\mg$.

It's also then straightforward to verify that the GIT quotient of the pluricanonical locus is a coarse moduli space for stable curves. An immediate corollary is that $\mgbar$ is projective. In fact, the construction depends both on the pluricanonical multiple $\nu$ used and on the sufficiently large auxiliary degree $m$ fixed in setting up the \GIT problem. Each pair of choices yields a natural ample class on $\mgbar$ and I have taken this opportunity to write down, in~\eqref{AmpleClasses}, formulae for these classes---they all lie in the $\lambda-\delta$-plane---that were worked out  in~\cite{HassettHyeonFlip} and~\cite{SwinarskiThesis} using ideas of Mumford~\cite{MumfordEnseignement} and Viehweg~\cite{ViehwegWeakI}.

\subsection{The Potential Stability Theorem}\label{potentialstabilitysubsection}

Fix $ g\geq 2 $ and a degree $d$. In the sequel, we will fix  a dimension $N$ implicitly determined by $d$ via Riemann-Roch as $N := d-g+1$.
Let $V$ be an $N$-dimensional vector space and let $ \h$ be the Hilbert scheme of curves in $\PV$ of degree $d$ and genus $g$ (or more precisely, with Hilbert polynomial $ \p(m)=md-g+1 $). Let $ \phi:X\to \h $  and $ L=\o_{H}(1) $ be the corresponding universal curve and  universal line bundle; we will abuse language and also write $X$
and $L$ for their restrictions to subschemes of $\h$.

\beginS{definition}%
\label{PotentiallyStable}%
 We call a connected curve $ C $ of genus~$ g $ and degree~$ d $ in $ \PV $ \emph{potentially stable} if: 
\beginS{enumerate} 
	\item The embedded curve $ C $ is nondegenerate (i.e., spans $\PV$).
	\item The abstract curve $ C $ is Deligne-Mumford semistable. 
	\item The linear series embedding $ C $ is complete and nonspecial: i.e., $ h^0(C, L)=N $ and $ h^1(C, L)=0 $. 
	\item \label{SubcurveInequality} [Subcurve Inequality] If $ Y $ is a complete subcurve of $ C $ of arithmetic genus~$ g_{Y} $ meeting the rest of $ C $ in $ k_{Y} $ points, then $$ \left|\deg_{Y}(L)-\frac{d}{g-1}\left(g_{Y}-1+\frac{k_{Y}}{2} \right)\right|\le\frac{k_{Y}}{2}. $$ 
\end{enumerate}
\end{definition}

The Subcurve Inequality implies that any chain of smooth rational components of $ C $ meeting the rest of $ C $ in exactly two points consists of a single smooth rational component embedded as a line. 
Thus, the abstract curve $ C $ underlying any potentially stable curve in $ \PV $ can only fail to be Deligne-Mumford stable in a very restricted way. We will continue, as usual, to abuse language and speak of a potentially stable curve $ C $ when the implied embedding is clear from the context. The justification for this somewhat baroque definition lies in the following theorem.

\beginSS{theorem}{[Potential Stability Theorem]} \label{PotentialStabilityTheorem} 
Fix integers $g$ and $d$ with $ g\geq 2 $. Suppose that $ d>9(g-1) $, or equivalently, that $\frac{d}{N} < \frac{8}{7}$. Then there is an $ M $ depending only on $ d $ and $ g $ such that if $ m\geq M $ and $ C $ in $ \PV $ is a connected curve with semistable $\thst{m}{th}$ Hilbert point, then $ C $ is potentially stable. 
\end{theorem}

The \pst\ came as a surprise when it first appeared since stable curves in the plane and other low-dimensional projective spaces can have arbitrarily bad singularities for large $g$. What Gieseker realized was that imposing the degree hypothesis above on the embedding does away with these pathologies. The proof of the \pst\ is the lengthiest step in the constructions of \mgbar\ and of the other moduli spaces we'll look at later, and can sometimes involve tedious technicalities. Despite the complications that ensue, the essential strategy is very simple: if $ C $ fails to have some property covered by Definition~\ref{PotentiallyStable}, find the filtration $ F $ of $ V $ that highlights this failure most clearly and check that $ F $ is destabilizing by showing some form of the numerical criterion  is violated. Only a certain care is needed in the order in which the properties are established since it is often necessary to assume some of these properties to justify estimates needed to verify that the failure of others is destabilizing. 

I'll lay out the sequence of steps here so that, when discussing newer constructions in later sections, I can focus on the points at which the arguments differ from Gieseker's model. I'll also sketch proofs of a few steps, likewise in preparation for discussions of the changes needed in other constructions. However, a complete proof is much too long to give here. For all the details, see \citelist{\cite{GiesekerTata} \cite{GiesekerCIME}} and for a somewhat condensed version \cite{Moduli}*{Section 4.C}. One definition is needed: if $Y$ and $Y'$ are any two subcurves of $C$ with no common components, let $k_{Y,Y'}$ denote the number of nodes at which $Y$ and $Y'$ intersect and let $k_Y= k_{Y, \overline{C\setminus Y}}\,$.

The proof of the Potential Stability Theorem proceeds via the following steps.
\beginS{enumerate}
	\item $\ired{C}$ is nondegenerate.
	\item Every component of $C$ is generically reduced.
	\item \label{rncstep} If an irreducible subcurve $Y$ of $C$  is 
	 not a rational normal curve, then $ \deg_Y\bigl(L\bigr) \geq 4 $.
	\item \label{cuspstep} If $Y$ is a reduced irreducible subcurve of $C$ then its normalization map $ Y_{ns} \to Y $ is unramified.
	\item Every singular point of $\ired{C}$ has multiplicity~2.
	\item \label{onlynodesstep} Every double point of \ired{C} is a node.
	\item $ H^{1}(\ired{C},L)=\{0\} $.
	\item $C$ is reduced, so $ H^{1}(C,L)=\{0\}$
	 and $V=H^{0}\bigl(C,L\bigl)$.
	\item For every subcurve $Y$ of $C$ and every 
	 component $E$ of the normalization $Y_{ns}$, either $ \deg{E}\geq k_{E,Y} $, or, $E$ is a rational normal curve for which $  \deg_E\bigl(L\bigr)=k_{E,Y}-1$.
	\item The Subcurve Inequality \ref{PotentiallyStable}.\ref{SubcurveInequality} holds for every subcurve $Y$ of $C$.
\end{enumerate}

As an illustration, let's look at Step~\ref{cuspstep}. First suppose that $p$ is an ordinary cusp. Consider the four-stage weighted filtration $F$ that gives weight~$0$ to the space $V_{3}$ of sections whose image under restriction to $Y$ and pullback via $\normaly$ to $Y_{ns}$ lie in $H^{0}\bigl(Y_{ns},L_{ns}(-4p)\bigr) $, weight~$1$ to the  space $ V_{2} $ of sections with images in $ H^{0}(Y_{ns}, L_{ns}(-3p)) $, weight~$2$ to the  space $ V_{2} $ of sections with images in $ H^{0}(Y_{ns}, L_{ns}(-2p)) $,  and weight~$4$ to all others. Since  $\normaly$ ramifies, $ Y $ itself must be singular.
Hence, $\deg_Y\bigl(\o_C(1)\bigr) \geq 4 $, by Step~\ref{rncstep}.

Then $\dim(V_0/V_1)=1$, $\dim(V_1/V_2) = 1$ and $\dim(V_2/V_3) = 1$ so the average weight $\alpha_F = \frac{7}{N}$. On the other hand, any $F$-monomial basis of $H^{0}\bigl(C,L^{\otimes m}\bigl)$ will restrict to a spanning set for $H^{0}\bigl(\ired{Y},L^{\otimes m}\bigl)$; since all the weights of $F$ are positive, $w_F(m)$ is at least the minimum weight of such a spanning set.

The key observation is that any monomial of weight at most $4m-i$ restricts on $\ired{Y}$ to a section vanishing to order at least $i$ at $p$. Further, since $\deg_Y\bigl(L\bigr) \geq 4$ and $Y$ is generically reduced, $\deg_{\ired{Y}}\bigl(L^{\otimes m}(-ip)\bigr) \geq 4m-i$. Thus Riemann-Roch gives the estimate
$$h^{0}\bigl(\ired{Y},L^{\otimes m}(-ip)\bigl) = m \deg_Y\bigl(L^{\otimes m}\bigr) -i - g + 1 + O(1)$$
with an implied constant depending only $g$. In other words, a spanning set must contain exactly one monomial of weight $i$ for almost all $i$ between $0$ and $4m$. Summing, we find that the weight of such a basis is at least $8m^2 + O(m)$ so has leading coefficent $e_F \ge 16$. 

Now recall from Lemma~\ref{ealphaInequality} of Section~\ref{AsymptoticNumericalCriterion} that $F$ is destabilizing if $e_F >  2 \alpha_F d$. Plugging in the values above, this is equivalent to $\frac{d}{N} < \frac{8}{7}$, exactly the numerical hypothesis of Theorem~\ref{PotentialStabilityTheorem}. 

If $p$ is a point of higher ramification, then we get a flag $F$ with $\alpha_F < \frac{7}{N}$, and the argument above again shows that $F$ is destabilizing.

Once we have Step 8 in hand, it's not hard to check that if $p$ is an ordinary cusp, then  $e_F = 16$. Hence $F$ is destabilizing if and only if $\frac{d}{N} < \frac{8}{7}$. Looking ahead to Section~\ref{schubert}, we note that $\nu$-canonically embedded curves have $\frac{d}{N} = \frac{2\nu}{(2\nu-1)}$ so, for these, $F$ is destabilizing if $\nu \ge 5$ but \emph{not} if $\nu < 5$. In fact, Hilbert stable curves \emph{can} carry ordinary cusps when $\nu \le 4$ as we'll see in Section~\ref{schubert}. 

The second case where the hypothesis $\frac{d}{N} < \frac{8}{7}$ is sharp arises in Step~\ref{onlynodesstep} when $C$ is the union of a curve $D$ of genus $(g-2)$ (not necessarily irreducible) embedded by a line bundle of degree $(d-1)$ and a tangent line $L$ to $D$ at a smooth point $p$. The singularity at $p$ is a tacnode but differs from other tacnodes in which neither local branch is a line in that the drop in degree on projecting from the tangent line $L$ is only $3$ rather than $4$.

In the latter case a flag giving weight $4$ to general sections, weight $2$ to sections vanishing at $p$ and weight $0$ to those vanishing on $L$ has, by similar arguments, $e_F \ge 16$ and $\alpha_F= \frac{6}{N}$ so is destabilizing if $\frac{d}{N} < \frac{4}{3}$. This suggests that bicanonically embedded curves with tacnodes can be Chow semistable, and calculations of Hassett and Hyeon~\cite{HassettHyeonFlip} also discussed in Section~\ref{schubert} confirm this.

\subsection{Properties of the pluricanonical locus}
\label{pluricanonicalproperties}

We'll say that a curve is nodal if it has at worst nodal singularities. Because deformations of nodal curves are either smooth or nodal, the subset $U$ of $\h$ parameterizing connected curves with at worst nodal singularities is open in \h. Since $ X \to U $ is a family of nodal curves, it has a  relative dualizing sheaf $ \omega=\omega_{X/U} $. 

Now we want to specialize by fixing a canonical multiple $ \nu \geq 3 $---to ensure that, by~\cite{Moduli}*{Exercise 3.10}, $\omega^{\otimes \nu}$ is very ample on all stable curves $C$---and set $d :=  \nu (2g-21) = \deg(\omeganu{C})$ and $ N = \homeganu{C} $. Since $\nu$ is fixed, we'll omit subscript $\nu$'s in denoting loci like $\h$ that implicitly depend on this choice.

\beginS{definition}
\label{ktwiddle}
We define the locus $J$ of $\nu$-canonically embedded stable curves to be the closed subscheme of $ U $ over which the sheaves L and $ \omega_{X/U}^{\otimes \nu} $ are equal. More formally,  $J$ is the subscheme defined by the \thst{(g-1)}{st} Fitting ideal of $ R^{1}\phi_{*}(\omega_{X/U}^{\otimes \nu}
\otimes L^{-1}) $.
\end{definition}

The locus $J$ is closed in $U$, hence locally closed in $\h$, and has dimension $ (3g-3)+(N^{2}-1) $: $(3g-3)$ for the choice of the underlying stable curve $C$ and $(N^{2}-1)$ for the choice of a basis of $\Homeganu{C}$ modulo scalars. Finally, $J$ is smooth: see~\cite{Moduli}*{Lemma 3.35}.

To go further we need to use Theorem~\ref{PotentialStabilityTheorem}. In terms of $\nu$, $\frac{d}{N} := \frac{2\nu}{(2\nu-1)}$ so the hypothesis $\frac{d}{N} < \frac{8}{7}$ of the theorem now requires $\nu \ge 5$. Thus, we know that every curve whose Hilbert point lies in the semistable locus $\semi{\h}$ of  $\h$ is potentially stable. The main claim is: 

\beginS{proposition} 
\label{ClosureOfKbar} $\semi{J}$ is closed in $\semi{\h}$.
\end{proposition}
\beginS{proof}	
We need an alternate form of the Subcurve Inequality \ref{PotentiallyStable}.\ref{SubcurveInequality} that follows by plugging in the formula
$\deg_{Y}(\restrictedto{\omega_{C}}{Y})=\deg_{Y}
 (\omega_{Y})+k_{Y}$ (or 
see \cite{Moduli}*{Exercise 4.47.2}):
\beginS{equation}
\label{SubcurveInequalityVariant}
 \left|\dy-\left(\frac{d}{\deg_{C}(\omega_{C})}\right)
 \deg_{Y}(\restrictedto{\omega_{C}}{Y})\right|\leq\frac{k_{Y}}{2}.
\end{equation}

Because $J$ is locally closed in the full Hilbert scheme $\h$, it is locally closed in  $\semi{\h}$. Applying the valuative criterion for properness, we must therefore show that given a discrete valuation ring $ R $ with residue field $k$ and quotient field $ F $, any map $ \alpha:\spec(R)\to\semi{\h} $ that takes the generic point $ \eta=\spec(F) $ of $ \spec(R) $ into $\semi{J}$ also takes the closed point $ 0=\spec(k) $ of $ \spec(R) $ into $\semi{J}$.
 
We first use $ \alpha $ to pull back the universal curve $ C \to\semi{\h} $ and the tautological bundle $L$ on it to $\spec(R) $. Let $ \omega=\omega_{C/\spec(R)} $ denote the relative dualizing sheaf of this family. It follows from the definition of $\semi{J}$ and the universal property of $\h$ that $ \alpha(0) $ will lie in $\semi{J}$ if and only if we can extend this isomorphism over the closed point $0$. 

The definition of $J$ also implies that $\restrictedto{L}{C_{\eta}}\cong\restrictedto{\omega^{\otimes \nu}}{C_{\eta}}$. Hence, if we decompose the special fiber $ C_{0} $ of $C$ into irreducible components $ C_{0}=\bigcup^{l}_{i=1}C_{i}$, then $ 
L \cong\omega^{\otimes \nu}\bigl(-\sum^{l}_{i=1}a_{i}C_{i}\bigr) 
$ 
with the multiplicities $ a_{i} $ determined up to a common integer translation. (Since $ \spec(R) $ is affine, $ \o_{C}(-C_0)\cong\o_{C} $.) We normalize the $ a_{i} $'s so that all are nonnegative and at least one equals~$0$. 

What we must show, then, is that \emph{all} the $ a_{i} $'s are 0. Note that this is automatic if $ C_{0} $ is irreducible. To take care of reducible $ C_{0} $'s, we use~\eqref{SubcurveInequalityVariant}. Let $ Y $ be the subcurve of $ C_{0} $ consisting of all $ C_{i} $ for which $ a_{i} $ is zero, and let $Z$ be the remainder of $ C_{0} $---i.e., those components for which $ a_{i} $ is positive. Then a local equation for $ \o_{C}\bigl(-\sum^{l}_{i=1}a_{i}C_{i}\bigr) $ is identically zero on every component of $Z$ and on no component of $ Y $. In particular, such an equation is zero at each of the $ k_{Y} $ points of $Y\cap Z$. Therefore, we find that
\vskip-18pt\beginS{equation*}
\beginS{split}
k_{Y} &\leq \deg_Y\Bigl(\o_C\bigl(-\sum^{l}_{i=1}a_{i}C_{i}\bigr)\Bigr)\\
		&=\deg_Y\bigl(\restrictedto{L}{C_0}\bigr)
 			-\nu\deg_Y\bigl(\restrictedto{\omega_{C/R}}{C_0}\bigr)\\
		&=\deg_Y\bigl(\restrictedto{L}{C_0}\bigr)
	-\Biggl(\frac{\deg_{C_0}\bigl(\restrictedto{L}{C_0}\bigr)}
	{\deg_{C_0}\bigl(\restrictedto{\omega_{C/R}}{C_0}\bigr)}\Biggr)
	\deg_Y\bigl(\restrictedto{\omega_{C/R}}{C_0}\bigr)\\
	&\le \frac{k_{Y}}{2}
\end{split}
\end{equation*}
where the last inequality follows from~\eqref{SubcurveInequalityVariant}.  Therefore $ k_{Y}=0 $ and since $ C_{0} $ is connected, $ a_{i}=0 $ for all~$ i $.
\end{proof}

\beginS{corollary} 
\label{StabilityAndJss} 
\beginS{enumerate}
	\item Every curve $ C $ in $\PV$ whose Hilbert point lies in $\semi{J}$ is Deligne-Mumford stable. 
	\item \label{semistablereductionstep} \semi{J} contains the $\nu$-canonical Hilbert point of every Deligne-Mumford stable curve of genus $g$. 
	\item $ \semi{J}=J^{s} $: every curve whose Hilbert point lies in $\semi{J}$ is Hilbert stable. 	
\end{enumerate}
\end{corollary}
\beginS{proof}
Every curve $ C $ in $\semi{J}$ is potentially stable so to prove (1) we only need to rule out smooth rational components meeting the rest of the curve in only two points. This is easy. On the one hand, the degree of the dualizing sheaf $ \omega_{C} $ of $ C $ on such a component is zero while, on the other, $ \omega_{C}^{\otimes n} $ is very ample on $ C $ because the Hilbert point of $ C $ lies in~\semi{J}.

For any Deligne-Mumford stable curve, $ \omega_{C}^{\otimes \nu}$ is very ample on $ C $, and thus embeds $ C $ as a curve in $\PV$ whose Hilbert point $ [C] $ lies in $\h$. To see that $ [C] $ lies in $\semi{J}$ or, equivalently, in $\semi{\h}$, choose a one-parameter deformation $ C\to\spec(R) $ of $ C $ to a smooth connected curve over a discrete valuation ring $ R $; that is, the generic fiber $ C_{\eta} $ of $C$ is a smooth curve of genus~$ g $ and the special fiber is $ C $. Then $C$ is again a stable curve over $ \spec(R) $, so its $ n $-canonical embedding realizes it as a family of curves in $ \PV $ over $ \spec(R) $ and hence corresponds to a unique morphism $ \alpha:\spec(R)\to\h $. Since the generic fiber $ C_{\eta} $ is smooth, its Hilbert point $ [C_{\eta}] $ lies in $\semi{\h}$ by Theorem~\ref{StabilityOfSmoothCurves}. This is the only, but essential, point at which this theorem is used in the whole construction.

We now obtain~\eqref{semistablereductionstep} by a \GIT semistable replacement argument (cf.~\cite{MumfordEnseignement}*{Lemma~5.3}). Since the quotient of $\semi{\h}$ by $\ssl(N) $ is projective, we can, after possibly making a finite change of base $\pi:\spec(R')\to \spec(R)$, find a map $ \beta:\spec(R')\to\semi{\h} $  such that the generic fiber $C'(\eta')$ of the pullback $C'$ of the universal curve over $\h$ by $\beta$ lies in the $\ssl(N)$-orbit of $ C_{\eta} $. By the uniqueness of the semistable reduction of a family of Deligne-Mumford stable curves, the stable models of the special fibers $ C_0 $ and $ C'_0 $ are isomorphic.  Since $ \beta(\eta') $ lies in $\semi{J}$, $ \beta(0') $ lies in $\semi{\h}$, and $\semi{J}$ is closed in $\semi{\h}$, we conclude that $ \beta(0') $ also lies in $\semi{J}$. In other words, $ C'_0 $ is also $\nu $-canonically embedded and hence must be Deligne-Mumford-\emph{stable}. Thus $ C_0 $ and $ C'_0 $ are both abstractly isomorphic and projectively equivalent in $ \PV $. But the Hilbert point $ [C'_0] $ is in $\semi{J}$ by construction, hence so is that of $ [C_0] $.

Every curve $ C $ whose Hilbert point lies in $\semi{J}$ is, by definition, Hilbert semistable. If the Hilbert point $[C]$ were not stable, then the closure of its $ \ssl(N) $-orbit would contain a semistable orbit with stabilizer of positive dimension. Since every curve whose Hilbert point lies in $\semi{J}$ is nondegenerate, this orbit would correspond to a curve $ C' $ with infinitely many automorphisms, and since $\semi{J}$ is closed in $\semi{\h}$, the Hilbert point of $ C' $ would lie in $\semi{J}$. This contradicts (1) and (3) follows.
\end{proof}

The upshot is that isomorphism classes of stable curves of genus $g$ correspond bijectively to \GIT\ stable $PGL(V)$-orbits in $\semi{J}$. It now follows by standard arguments from the universal property of the Hilbert scheme that defining $\mgbar := \semi{J}\quotient PGL(V)$ is a coarse moduli space for such curves. See~\cite{NewsteadModuli}*{Proposition 2.13} for the general argument or~\cite{Moduli}*{p. 222} for this case.

\subsection{A few applications}
\label{applications}

This section is a quick review of some corollary information about $\mgbar$ that can be deduced from the \GIT\ construction. In this volume, the first consequence to mention must surely be the irreducibility of $\mgbar$ in positive characteristics because this, as is clear from its title, was the result that motivated the great paper of Deligne and Mumford~\cite{DeligneMumford}. The argument, based on irreducibility in characteristic $0$, can be found in~\cite{GiesekerTata}*{Lemma 2.0.3} or~\cite{Moduli}*{Theorem 4.43}}.

Next, the projectivity of $\mgbar$, first proved by Knudsen in the series of paper~\cites{KnudsenMumford, KnudsenII, KnudsenIII}, is an immediate consequence of the its construction as a \GIT\ quotient the closed subscheme $\semi{J}$ of $\semi{\h}$. Indeed, it comes equipped with a polarization, given by the $\proj$ of the ring of invariant sections of powers of the $\pgl(V)$-linearized bundle $\Lambda_m$ used on $\h$. Or, better, with many polarizations, since $\h$ depends on the choice of the canonical multiple $\nu$ and then $\Lambda_m$ depends on the choice of the sufficiently large degree $m$. These have recently been computed by Hassett and Hyeon~\cite{HassettHyeonFlip} and Swinarski~\cite{SwinarskiThesis} following Mumford's calculation in~\cite{MumfordEnseignement} for the $\nu$-canonical Chow quotient. We'll need them in Sections~\ref{schubert} and~\ref{open}.

I'll use additive notation for the operation in $\pic$ for legibility. To start with denote by $C \subset \semi{J} \times \PV$ the universal curve over  \semi{J} and let $\pi$ be the projection of $C$ onto the first factor. Since the fibers of $C\to \semi{J}$ are $\nu$-canonically embedded, we must have 
\beginS{equation}
\label{definingrelationforQ}
\o_{C}(1) = \nu \cdot \omega_{C/\semi{J}} + \pi^*(Q)
\end{equation}
for some invertible sheaf $Q$ on $\semi{J}$.  To solve for $Q$, we can take the direct image of~\eqref{definingrelationforQ} and apply the fact that  $(\pi_*\bigl(\o_{C}(1)\bigr)$ is a trivial bundle of rank $N$. This leads to  
\beginS{equation}
\label{solveforQ}
- N \cdot Q = c_1\Bigl(\pi_*\bigl(\nu \cdot \omega_{C/\semi{J}}\bigr)\Bigr) = \binom{\nu}{2} \kappa +  \lambda\, , 
\end{equation}
where the last equality follows from the  Grothendieck-Riemann-Roch calculation of \cite{MumfordEnseignement}*{Theorem 5.10}. 

Taking the $\thst{m}{th}$ tensor power of~\eqref{definingrelationforQ} using the right inequality in~\eqref{solveforQ} a second time with $\nu$ replaced by $m\nu$ yields the first equality in
\beginS{equation}
\label{coneoCm}
\beginS{split}
	\Lambda_m
	&= \binom{m\nu}{2} \kappa + \lambda + \frac{m P(m)}{N}Q \\
	&= \binom{m\nu}{2} \kappa + \lambda - \frac{m (2\nu m-1)}{(2\nu-1)} \bigl(\binom{\nu}{2} \kappa +\lambda\bigr)
\end{split}
\end{equation}
and the second follows by substituting for $Q$ from~\eqref{solveforQ} and writing $P(m)$ and $N$ in terms of $\nu$ and $m$. 

This is the polarization that our construction for $\nu$ and $m$ yields and so gives an ample class on $\mgbar$ whenever our construction can use these values. Note that it is already evident that \emph{none} of these classes depends on $g$. Applying the relation $\kappa = 12 \lambda - \delta$ and doing some algebraic simplification we find:
\beginS{lemma}\label{AmpleClasses}
For $\nu \ge 2$ and $m$ large enough that the quotient of the set of $\thst{m}{th}$-Hilbert points of $\nu$-canonically embedded stable curves is $\mgbar$, the induced very ample polarization has class
\beginS{equation}\label{polarizationformula}
\left(\frac{m-1}{2 \nu -1}\right) \left(\bigl( 6 \nu^{2}m-2\nu m-2\nu+1 \bigr) \lambda
-\frac{\nu^2 m}{2}\delta\right).
\end{equation}
\end{lemma}
The slope of these divisors, by which I mean the ratio of the $\lambda$ and $-\delta$ coefficients, is therefore $12 - \frac{4}{v} - \frac{4}{vm}+ \frac{2}{v^2m}$. For fixed $\nu$, these values approach $12- \frac{4}{v}$ from below as $m \to \infty$. As we vary both $m$ and $\nu$, they range over an interval $[11.2-\epsilon, 12)$: the value of epsilon depends on how large we need to take $m$ when $\nu = 5$ and the range does \emph{not} include $12$. 

As a check, we recall that a theorem of Mumford and Knudsen~\cite{KnudsenMumford}*{Theorem 4} (see also~\cite{FogartyTruncated}), that I'll write only for our situation, says that if $\chow: \semi{J} \to \ddiv$ is the Chow map, then there are invertible sheaves $\mu_0$, $\mu_1$ and $\mu_2$ on $\semi{J}$ such that, $\mu_2 = \chow^*\bigl(\o_{\ddiv}(1)\bigr)$ and, for large enough $m$, 
\beginS{equation}
\label{knudsenequation}
\Lambda_m = \sum_{i=0}^2 \binom{m}{i} \mu_i
\end{equation}
We can for solve $\mu_2$ by equating $m^2$-coefficients in this equation obtaining $\mu_2 = \frac{\nu}{2\nu-1}\bigl((12\nu-4)\lambda - \nu\delta\bigr)$ which recovers the polarizations computed by Mumford using Chow quotients in~\cite{MumfordEnseignement}*{Corollary 5.18}. Plugging in and repeating, we also find that $\mu_1=-\lambda$ and that $\mu_0=\lambda$.

In the other direction, a class $a\lambda - b\delta$ cannot be ample if the slope $s = \frac{a}{b}\le 11$. Recall that elliptic tail is a \emph{connected} genus $1$ subcurve of a stable curve meeting the rest of the curve at a single point. A varying elliptic tail is a curve $\PP^1 \subset \mgbar$ obtained by gluing an elliptic tail of varying $j$-invariant to a fixed point on a fixed curve of genus $g-1$. An explicit family can be obtained, for example, by blowing up a generic pencil of plane cubics at the $9$ base points and taking one of the exceptional divisors as the section determining the marked point. It's then straightforward to check that $11\lambda - \delta$ has degree $0$ on such a curve (see~\cite{Moduli}*{Example 3.140}). Therefore in $\mgbar$, this class contracts the divisor $\Delta_1$. The question of where in the gap between slope $11.2$ and $11$ the boundary of the intersection of the ample cone of $\mgbar$ with the $\lambda-\delta$ plane lies was settled by Cornalba and Harris who showed that $s >11$ is also sufficient. Their proof is based on a positivity result that we mention here because it depends on the ideas of Section~\ref{numerical}.

\beginSS{theorem}{[Cornalba-Harris Theorem~\cite{CornalbaHarrisAmple}]}
\label{cornalbaharristheorem}
Suppose  $\pi: X \to B$ is a proper flat family of relative dimension $r$ over an irreducible curve $B$ and $L$ is a line bundle on  $X$ such that $\pi_{*}L$ is a vector bundle $E$ of rank $N$. Suppose further  that, for a general point $ b\in B$, the line bundle $\restrictedto{L}{X_{b}}$ is very ample and embeds $X_{b}$ as a Hilbert stable variety in $\PP^{N-1}$.  Then 
$$N\cdot \pi_*\bigl(c_{1}(L)^{r+1}\bigr)\geq(r+1)\cdot \pi_*\bigl(c_{1}(L)^{r}\bigr)\cdot c_{1}(E)\, .$$
\end{theorem}

The requirement of stability in this theorem seems, at first glance, strange. In the proof, it guarantees the existence of a homogeneous invariant that is interpreted as a section of an auxiliary line bundle $M$ on $B$ and the inequality is deduced from the non-negativity of the degree of $M$. Moreover, an example due to the author shows that without this hypothesis the inequality of the Theorem may fail. The main step in pinning down the slope of the ample cone is to apply the inequality when $X\to B$ is a family of curves with smooth, non-hyperelliptic general fiber to deduce 
\begin{equation}\label{smoothnonheinequality}
	\deg_B(\lambda) \ge \Bigl(8 + \frac{4}{g}\Bigr) \deg_B(\delta)\,.
\end{equation}
Families of curves whose general member is hyperelliptic and/or singular were handled by arguments that do not involve stability (for which see~\cite{Moduli}*{6.D}). Recently, Lidia Stoppino~\cite{StoppinoSlopeInequalities} proved a variant of the Cornalba-Harris theorem that applies to families whose general fiber is smooth and hyperelliptic.

The Hilbert stability hypothesis needed to apply Theorem~\ref{cornalbaharristheorem} to get~\eqref{smoothnonheinequality} is provided by,
\beginS{lemma}
\label{StabilityOfCanonicalNonHECurves} 
If $C$ is a smooth, non-hyperelliptic curve of genus $g\geq 2$ embedded in $\PP^{g-1}$ by its canonical linear series $K$, then $C$ is asymptotically Hilbert stable.	
\end{lemma}
\begin{proof} The proof is almost identical to that for Theorem~\ref{StabilityOfSmoothCurves} where the only geometric ingredient needed was the inequality $\epsilon_{i} < \big(\frac{d}{N-1}\bigr)i$. For the canonical series, this becomes $\epsilon_{i} < \bigl(\frac{2g-2}{g-1}\bigr)i = 2i$. Referring back to Figure~\ref{RiemannRochLineDiagram}, we thus need to rule out the existence of sub-linear series on the Clifford line and these exist if and only if $C$ is hyperelliptic. 
\end{proof}

In Stoppino's argument, a theorem of Kempf's~\cite{KempfInstability}*{Corollory~5.3} which implies that any rational normal curve has semistable Hilbert point substitutes for this result.


\section{Pseudostable Curves}
\label{schubert}

\subsection{Constructions of $\nu$ canonical quotients for $\nu <5$}\label{smallnuoverview} 

In the preceding section, we exhibited $\mgbar$ as the quotient of the locus $J$ in the Hilbert scheme $\h$ of $\nu$-canonical curves. Since, for a Deligne-Mumford stable $X$, $\omega_X$ is very ample when $\nu \ge 3$, a natural question is to describe the quotient when $\nu$ equals $3$ or $4$.

The hypothesis $\nu \ge 5$ is sharply invoked in Gieseker's construction only to apply the calculation that ordinary cusps destabilize Hilbert and Chow points of curves of degree $d$ in $\PP^{N-1}$ when $\frac{d}{N}< \frac{9}{8}$. See the discussion of Step~\ref{cuspstep} in the proof of the Potential Stability Theorem (Theorem~\ref{PotentialStabilityTheorem})---the tacnodal curves arising in Step~\ref{onlynodesstep} for which this inequality is also sharp do not arise as pluricanonical models. The calculations there suggest (but do not prove) that for $\nu < 5$ some cuspidal curves will be stable. Via this the hypothesis $\nu\ge5$ enters implicitly into the proof that nodal stable curves are Hilbert stable. 

Another clue is provided by comparing the Cornalba-Harris Theorem~\ref{cornalbaharristheorem} with the computations of polarizations in Lemma~\ref{AmpleClasses}. Note that, for $\nu =4$, the slope approaches $11$ from below as $m \to \infty$; the corresponding values are $\frac{32}{3}$ for $\nu=3$, $10$ for $\nu=2$,  and $8$ for $\nu=1$ and these limits give the slopes of the polarizations on the corresponding Chow quotients. In particular, we can already see that for $\nu <5$, none of these \GIT\ quotients can be $\mgbar$ since, at the least, the divisor $\Delta_1$ must be contracted. 

These observations suggest that ordinary cusps appear and elliptic tails disappear in the $3$- and $4$-canonical quotients.
This prediction was verified by Schubert~\cite{Schubert} who, using $3$-canonical Chow points, produced a quotient that is a coarse moduli space $\mpsgbar$ for pseudostable curves for $g \ge 3$.

\beginS{definition}\label{pseudostabledefinition}
A curve is \emph{pseudostable} if:
\beginS{enumerate}
		\item It is reduced, connected and complete.
		\item It has finite automorphism group.
		\item Its only singularities are nodes and ordinary cusps.
		\item It has no elliptic tails.
\end{enumerate}
The requirement of having a finite automorphism group means that any  component whose normalization is rational must contain at least $3$ singular points. If a curve meets the other requirements but has components with rational normalization containing $2$ singular points, then it is \emph{semipseudostable}. In addition to the chains of rational curves that Deligne-Mumford semistability permits, rational cuspidal tails (rational cuspidal curves meeting the rest of the curve in a single point) are semipseudostable.
\end{definition}

\subsection{Applications to the log minimal model program} 
\label{logminimal}
Before I discuss Schubert's proof, I want to mention other recent work that involves understanding $\nu$-canonical quotients for even smaller $\nu$. Most of this arises in connection with the log minimal model program for $\mgbar$ (and recently other spaces $\mgnbar$) initiated by Hassett and Hyeon. Recall that the goal here is to understand the model $\mgbar(\alpha)$ of $\mgbar$ that arises as $\proj\bigl(\oplus_{n\ge 0} \Gamma(n(K_{\Mgbar}+\alpha\delta))\bigr)$. Most of the results in this area to date proceed indirectly. First GIT is used to construct a quotient that carries the natural polarization proportional to $K_{\Mgbar}+\alpha\delta$ and then this quotient is identified with $\mgbar(\alpha)$. Although I only give details of Schubert's construction below because it is the simplest model, I'd like to point out some interesting new ideas in the more recent constructions.

Gieseker's construction can be viewed as the implementation of the log minimal model program for $1 \ge \alpha > \frac{9}{11}$. In the paper~\cite{HassettHyeonLogCanonical}, Hassett and Hyeon show that Schubert's $3$-canonical quotient is $\mgbar(\alpha)$ for $\frac{9}{11} \ge \alpha > \frac{7}{10}$. More recently in~\cite{HassettHyeonFlip}, they have constructed, for $g\ge 4$,  $\mgbar(\frac{7}{10})$ as the quotient of the Chow variety of bicanonical curves and $\mgbar(\frac{7}{10}-\epsilon)$ as the quotient of the Hilbert scheme of bicanonical curves. 

This quotient turns out to be a moduli space $\mhgbar$ for a class of curves that they dub h-semistable in which nodes, cusps and tacnodes are allowed but certain chains of elliptic curves are excluded (see~\cite{HassettHyeonFlip}*{Definitions~2.4-2.6} for details). However, the instability calculations are quite a bit trickier than Schubert's and, for most $g$, there are strictly semistable points.  
\setdiagram{flipdiagram}{nohug,grid=flipgrid}{
\mgbar({\scriptstyle\frac{7}{10}}+\epsilon)& & & & \mgbar({\scriptstyle\frac{7}{10}}-\epsilon)\\
& \rdTo^{\Psi}& &\ldTo^{\Psi^+} & \\
& & \mgbar({\scriptstyle\frac{7}{10}})& & \\
}%
Moreover these spaces fit into a picture like that shown in~\eqref{flipdiagram} in which $\Psi$ is a small contraction and $\Psi^+$ is its flip.

In~\cite{HyeonLeeGenusThree}, Hyeon and Lee complete the log minimal model program for $g=3$ producing an analogous flip at $\alpha=\frac{7}{10}$, and then using a GIT analysis of plane quartics (the canonical models in this case) to show that the only other critical values are $\frac{17}{28}$---here the hyperelliptic locus gets contracted---and $\frac{5}{9}$ at which the whole space is contracted to a point.

I'd like to mention one novel stability criterion in ~\cite{HassettHyeonFlip}*{section 3}. The idea is that the relation~\eqref{knudsenequation} for the linearization $\Lambda_m$ in terms of the tautological classes $\mu_i$ (denoted $L_i$ in~\cite{HassettHyeonFlip}) holds for \emph{small} $m$ on the locus of Hilbert points of curves for which restriction of homogeneous polynomials of degree $m$ is onto $\HO{C}{\o_C(m)}$ and for which $\o_C(m)$ has no higher cohomology--in particular, where $C$ is $m$-regular in the sense of Castelnuovo-Mumford. On this locus, we therefore get the same relation between the least $\rho$-weights of $\Lambda_m$ and the $\mu_i$ for any $1$-ps $\rho$. If moreover, $\o_C(1)$ has no higher cohomology and $V=\HO{C}{\o_C(1)}$ is an isomorphism, then $\mu_0=-\mu_1$ \cite{HassettHyeonFlip}*{Proposition 3.9}. This forces the divisibility by $(m-1)$ seen in~\eqref{polarizationformula}.

For such a curve, the least $\rho$ weight in any 2 degrees greater than or equal to the regularity determine the polynomial give the $\Lambda_m$-weight for all large $m$. In particular, if $C$ is $2$-regular, Proposition~3.17 shows that 
\begin{equation}\label{mtwothree}
w_{\rho}(m) = (m-1)\bigl((3-m)w_{\rho}(2)+ (\frac{m}{2}-1)w_{\rho}(3)\bigr)
\end{equation}
In particular, this formula makes is possible to use tools like Macaulay~2 to automate many instability checks: see also~\cite{HassettHyeonLee} for other examples.

Finally, I want to mention the beautiful paper~\cite{HassettGenusTwo} of Hassett in which he deals with the case $g=2$ by techniques which are special to that case and which inaugurated work in this area. In particular, he uses the explicit invariant theory of binary sextics to describe the various log minimal models, realizing, for example, the model $\mibar{2}(\frac{9}{11})$ as the resulting projective quotient.

\subsection{Overview of Schubert's Proof}
\label{schubertoverview}

Schubert's argument follows the general lines of Gieseker's. I'll sketch it briefly here highlighting the points of significant difference and then return to discuss these in more detail later in this section. 

First come various stability results.  Schubert needs to know that smooth curves have stable Chow points which again follows from Theorem~\ref{StabilityOfSmoothCurves}. He also proves a pseudostable variant of the Potential Stability Theorem for $3$-canonical models that differs from Theorem~\ref{PotentialStabilityTheorem} in two ways matching the expectations above. Both nodes and ordinary cusps are now allowed and elliptic tails are shown to be destabilizing, hence are prohibited. I'll return to this last point in a moment.

The major novelties in the argument are substitutes for standard theorems about pairs of families of stable curves over a discrete valuation ring having isomorphic smooth generic fibers. The first (his Lemma~4.2, here~\eqref{separatednessofpseudostable}) can be viewed as a valuative criterion of separatedness for the functor of flat families of pseudostable curves. It asserts that if the special fibers of both families are pseudostable, then they are isomorphic. The second (Lemma~4.8, here~\eqref{dmversuspseudostable}) functions as a substitute for semi-stable reduction in Corollary~\ref{StabilityAndJss}.\ref{semistablereductionstep} and is the key to proving that $3$-canonical pseudostable curves have stable Chow points. It asserts that if one family has a Deligne-Mumford stable special fiber and the other has a pseudostable special fiber, then there is a map from the stable family to the pseudostable one that is an isomorphism except over cusps of the pseudostable special fiber above which an elliptic tail is contracted. I'll sketch the ideas behind these results in~\ref{schubertlemmas}.
	 
Dave Swinarski pointed out to me that Schubert says nothing about  the case $\nu=4$, whose GIT you'd expect to be easier, and the paper~\cite{HyeonMorrison} explains why. It turns out that Schubert's argument for $3$-canonical Chow points applies with only obvious adjustments to $3$-canonical Hilbert points and that most of it applies also to $4$-canonical Chow and Hilbert points. Only his Lemma~3.1 showing that curves with elliptic tails are \GIT\ unstable breaks down: his argument applied to $4$-canonical curves only shows that such curves are not Chow stable and says nothing about their Hilbert stability. However, by specifying the $1$-parameter subgroup $\lambda$ used in his argument a bit more carefully and making a more precise analysis of the weights with $\lambda$ acts it is possible to show (see Corollary~\ref{fourcanonicaltailsunstable}) that it destabilizes $4$-canonical Hilbert points. In~(\ref{elliptictails}), I review the argument.

Substituting this result for his Lemma~3.1, the remainder of Schubert's construction, with the Chow scheme for $3$-canonical curves replaced by the  $\thst{m}{th}$ Hilbert scheme of either $3$- or $4$-canonical curves for a sufficiently large $m$, goes through with only minor changes. In view of this and for consistency with the rest of this paper, I'll stick to the Hilbert schemes version of Schubert's construction in what follows.  This leaves open the question of the $4$-canonical Chow quotient. Here the geometry is more complicated. There are three classes of strictly stable orbits that are identified in the quotient: both curves with elliptic tails and curves with cusps are in the basin of attraction of curves with a rational cuspidal tail. For more details, see~\cite{HyeonMorrison}.

\subsection{Pluricanonical Stability of Elliptic Tails}
\label{elliptictails}

This subsection gives the refinement of Schubert's analysis of stability of elliptic tails needed to make his construction apply to $4$-canonical models. The arguments follows closely that in~\cite{HyeonMorrison}. First, we recall the setup.

Fix a Deligne-Mumford stable curve $X$ with an elliptic tail, i.e. $X = C \cup E$ where $C$ and $E$ are subcurves of genera $(g-1)$ and $1$ respectively and $C\cap E$ is a single node $p$. Note that $C$ is not assumed to be smooth or irreducible. Assume $\nu \ge 3$ so that $\lbpow{\omega_X}{\nu}$ is very ample and let $d=2\nu(g-1)=\deg\big(\lbpow{\omega_X}{\nu}\bigr)$ and $N = d-g+1= \h^0(X,\lbpow{\omega_X}{\nu})$. Then $\lbpow{\omega_X}{\nu}$ has restriction to $E$ linearly equivalent to $\o_E(\nu p)$ and has degree $c=(d-\nu)-g+2$. To simplify notation, I'll write $L=\lbpow{\omega_X}{\nu}$.

It follows directly from Riemann-Roch that the linear spans $V_C$ of $C$ and $V_E$ of $E$ in $\PP^{N-1}$ are of dimensions $c-g+1 = N-\nu$ and $\nu-1$ respectively and that their intersection is $\{p\}$. Letting $l=N-\nu+1$, we can therefore choose homogeneous coordinates $x_1, \ldots, x_N$ such that $ x_1=  \ldots = x_{l-1} = 0$ defines $V_E$, $ x_{l+1}= \ldots = x_{N} = 0 $ defines $V_C$, and $p$ is the point where all the $x_i$ except $x_l$ vanish. 

For $ j \ge 1 $, we will confound $x_{l+j}$ with the section of $\HO{E}{\restrictedto{L}{E}}$ it determines and write $\ord_p(x_{l+j})$ for the order of vanishing at $p$ of this section. Again by Riemann-Roch, we may choose $x_{l+j}$ so that  $\ord_p(x_{l+j})=j$ for $1 \le j \le \nu-2$ and choose $x_N$ so that $\ord_p(x_N)=\nu$.

Define $\lambda$ to be the $1$-ps subgroup of $\ssl(N)$ acting by $\diag(t^{r_1}, \cdots, t^{r_N})$ in these coordinates where
$r_i$ equals $\nu$ if $i \le l$, $\nu -j$ if $i= l+j$ and $1 \le j \le \nu -2$ and $0$ if $i=N=l+\nu -1$. Note that, for $ j \ge 0$, this gives $x_{l+j}$ weight equal to $e-\ord_p(x_{l+j})$. 

The proof of Schubert's Lemma~3.1 shows that the $m^2$ coefficient of $w_{\lambda}(m)$ is at least $\bigl(d-\frac{\nu}{2}\bigr)\nu$ (although only $L=\lbpow{\omega_X}{3}$ is considered). The next lemma is the sharpening of this estimate to an exact evaluation of $w_{\lambda}(m)$ needed to apply his argument to $4$-canonical models.
\beginS{lemma}\label{lambdaweightlemma}
~~$\displaystyle{w_{\lambda}(m) = m^2\bigg[\bigl(d-\frac{\nu }{2}\bigr)\nu \bigg] + m\bigg[\bigl(\frac{3}{2}-g\bigr)\nu \bigg] -1}$.
\end{lemma}
\beginS{proof}
For concision, we will henceforth understand all monomials to have degree $m$ and view them directly as sections of $\lbm{L}$ over $X$ or $E$ (eliding ``the restriction to''). Weights will always be $\lambda$-weights. 

Let $W_r$ be the span in $\HO{X}{\lbm{L}}$ of all monomials of weight at most $r$ and let $s=m\nu -r$.   I claim that
\beginS{equation*}
	\dim(W_r)
	=\beginS{cases}
	md-g+1 & \text{if $r =  m\nu $}\\
	r & \text{if $2 \le r \le  m\nu -1$}\\
	1 & \text{if $r=0$ or $r=1$}
	\end{cases}		
\end{equation*}
Given this, the lemma follows by elementary manipulations since the weight of any monomial basis is simply the sum of $r \big(\dim(W_r)-\dim(W_{r-1})\bigr)$ over $r$.

The first case in the claim is immediate from Riemann-Roch for $\lbm{L}$ on $X$. The others follow from the equality
\beginS{equation}\label{weightspaces}
	W_{r} = \HO{E}{\restrictedto{\lbm{L}}{E}(-sp)} \text{~~for $r=0$ \text{and for} $2\le r \le m\nu -1$\,.}
\end{equation}
Since $\restrictedto{\lbm{L}}{E}(-sp)\cong \o_E(rp)$---recall that $\restrictedto{L}{E} \cong \o_E(ep)$, Riemann-Roch on $E$ implies that $\hO{E}{\restrictedto{\lbm{L}}{E}(-sp)}=r$.

If any monomial has weight $r=m\nu -s$ then it contains one or more factors $x_{l+j}$ with $j>0$ and hence vanishes on $C$. By construction, $s$ equals the sum of the orders of vanishing at $p$ of the factors of this type, hence $W_{r} \subset \HO{E}{\restrictedto{\lbm{L}}{E}(-sp)}$. 

If we next set $M_0 = x_n^m$, then $B_0 := \{M_0\}$ is a basis of $\HO{E}{\restrictedto{\lbm{L}}{E}(-m\nu p)}$ lying in $W_0$. Finally, for $r=2, \ldots, m\nu -1$, let $M_r$  be any monomial 
$$M_r := \prod_{k=1}^m x_{l+j_k} \text{~s.t. each $j_k \ge 0$ and} \sum_{k=1}^m j_k = s=m\nu -r\,.$$
Then, $M_r$ vanishes on $C$ because  some $j_k>0$.  By construction, $M_r$ has weight exactly $r$ and, since $x_l$ is non-zero at $p$, $M_r$ vanishes to order exactly $s$ at $p$. Thus, $B_r:= \{M_0, M_2, M_3, \ldots, M_r\}$ is a subset of $W_r \cap \HO{E}{\restrictedto{\lbm{L}}{E}(-sp)}$ of cardinality equal to $\hO{E}{\restrictedto{\lbm{L}}{E}(-sp)}$. But all the elements of $B_r$ except $M_r$ lie in $W_{r-1}$. By induction, $B_r$ is linearly independent and hence is a basis of $\HO{E}{\restrictedto{\lbm{L}}{E}(-sp)}$ which therefore lies in~$W_{r}$.
\end{proof}

We now want to apply the Numerical Criterion (\ref{NumericalCriterionForHilbertPoints}). In our examples, $P(m)=md-g+1$ by Riemann-Roch and an easy calculation shows that  
\beginS{equation*}
	\alpha(\lambda) = \frac{\nu N-\sum_{j=1}^{\nu -2} j - \nu }{N} = \nu - \frac{\nu ^2-\nu +2}{2N}
\end{equation*}
so we want to compare $w_{\lambda}(m)$ to
\beginS{equation}\label{mpmalpha}
m P(m) \alpha(\lambda) = m^2\bigg[d\bigl( \nu - \frac{\nu ^2-\nu +2}{2N}\bigr)\bigg]+	m\bigg[(1-g)\bigl( \nu - \frac{\nu ^2-\nu +2}{2N}\bigr)\bigg]
\end{equation}

After some easy simplifications, we find that the $m^2$ coefficient in Lemma~\ref{lambdaweightlemma} is less than that in~\eqref{mpmalpha} equal, if and only if    
$$\frac{d}{N} = \frac{2\nu}{(2\nu-1)} < \frac{\nu ^2}{(\nu ^2-\nu +2)}$$
and that this happens exactly when $\nu \ge 5$. In these cases, the \thst{m}{th}-Hilbert point of $X$ is $\lambda$-stable for large enough $m$. We get the opposite comparison (and hence an unstable Hilbert point) exactly when $\nu \le 3$. But for $\nu=4$ the two coefficients are equal. In this case, to decide the whether $X$ is Hilbert stable or unstable with respect to $\lambda$, we need only compare the two $m$ coefficients. Plugging in $\nu=4$ and simplifying, we find that 
$
(1-g)\bigl( \nu - \frac{\nu ^2-\nu +2}{2N}\bigr)-\bigl(\frac{3}{2}-g\bigr)\nu = -1
$
so the $4$-canonical Hilbert point of $X$ is unstable. Note that, as predicted by~\eqref{mtwothree}, $w_F(m)$ is divisible by---in fact, equal to---~$m-1$.

Thus we have proved the first two claims below.  I emphasize that, although we have shown that, when $\nu=4$, $X$ is Chow strictly stable with respect to $\lambda$, $X$ might be unstable with respect to some other $1$-ps. Likewise the third claim follows not from the preceding argument, showing that $X$ is $\lambda$-stable for $\nu \ge 5$, but from Corollary~\ref{StabilityAndJss}. 
\beginS{corollary}\label{fourcanonicaltailsunstable}
Let $X$ be a $\nu$-canonically embedded Deligne-Mumford stable curve with an elliptic tail. Then,
\beginS{enumerate}
	\item If $\nu =3$, $X$ is Chow unstable and asymptotically Hilbert unstable.	
	\item If $\nu =4$, $X$ is not Chow stable and is asymptotically Hilbert unstable.	
	\item If $\nu \ge 5$, $X$ is Chow stable and asymptotically Hilbert stable.	
\end{enumerate}
\end{corollary}	

\subsection{Schubert's Key Lemmas}
\label{schubertlemmas}

Here I want to sketch the proofs of the two lemmas about pseudostable curves cited in subsection~\ref{schubertoverview} and explain how these are applied in constructing $\mpsgbar$. First, to fix notation, let $R$ be a discrete valuation ring, let $B=\spec(R)$, and let $\eta$ and $0$ be the generic and special points of $B$. 

\beginS{lemma} \label{separatednessofpseudostable}
	If $\pi: Y\to B$ and $\pi':Y'\to B$ are flat families of pseudostable curves with isomorphic smooth generic fibers, then the special fibers $Y_0$ and $Y'_0$ are also isomorphic.
\end{lemma}
\beginS{proof}
The key claim is that: if $g \ge 3$, a stable curve $Z$ has a unique \emph{connected} subcurve $C$ containing no elliptic tails and whose complement consists of a set of pairwise disjoint elliptic tails.

Given this, the first step is a standard application of Stable Reduction to see that after a base change, if necessary, there is a family $Z\to B$ with generic fiber isomorphic to those of $Y$ and with stable special fiber $Z_0$ and a $B$-maps $\phi:Z\to Y$. As in the discussion in~\cite{Moduli}*{pp. 122-130}, $Z_0$ will have an elliptic tail of $j$-invariant $0$ lying over each cusp of $Y_0$. By pseudostability, these will be the only elliptic tails in $Z_0$; further, as there are no other non-nodal singularities and no unstable rational components, $Z_0$ will be isomorphic to $Y_0$ except over the cusps. More precisely, the complement of the points of attachment of the elliptic tails in $C_0$ will be isomorphic to the complement of the cusps in $X_0$ and there there will be one elliptic tail over each cusp.  

Now apply the same argument to $Y'$. The uniqueness of stable reductions implies that we get the same stable central fiber $Z_0$ and hence that $Y_0$ and $Y'_0$ are isomorphic away from their sets of cusps which are in canonical bijection. Hence $Y_0$ and $Y'_0$ are isomorphic. 

The claim is not quite as trivial as it may appear. We define $C$ inductively starting with $C=X$. If $C$ contains an elliptic tail $E$ not meeting any elliptic tail in the complement of $C$ then, replace $C$ with the closure of the complement of $E$. If not, stop. If $C$ is a connected genus $1$ subcurve of $X$ meeting the deleted elliptic tails in a single point then $X$ is the union two elliptic tails and has genus $2$. This is the first point where Schubert needs to assume $g\ge 3$: if $X$ is a general point of $\Delta_i$ in genus $2$ (the join of $2$ elliptic tails), there is no canonical subcurve $C$.  If not, any connected genus $1$ subcurve $E$ must meet both the rest of $C$, by connectedness, and the set of deleted elliptic tails, by induction, so is not an elliptic tail.

Next, I leave the reader to check that any elliptic tail is irreducible. Given this, if $C$ and $C'$ both satisfy the claim and $E$ is an elliptic tail deleted from $X$ in forming $C$, $E$ either lies inside $C'$, contradicting its choice, or lies in its complement. This shows that $C' \subset C$ and, by symmetry, proves uniqueness.
\end{proof}

To complete the construction we will need a Corollary of the proof above that is the content of Schubert's Lemma~4.2. Say that a stable curve $Z$ has \emph{standard tails} if all its elliptic tails are smooth with $j$-invariant $0$. Let $C$ be the canonical subcurve obtained as in the claim by deleting all these elliptic tails and let $D_Z= C \cap \overline{Z\setminus C}$. For a pseudostable curve $Y$, let $D_{Y}$ be the set if cusps on $Y$.
\beginS{corollary}\label{stablemodelofpseudo}
If $Y$ has genus $g\ge 3$, then there is a unique stable curve $Z$ with standard tails and a map $\pi:Z\to Y$ such that $\restrictedto{\pi}{C}$ is the normalization of $Y$, $\restrictedto{\pi}{C\setminus D_Z}$ is an isomorphism to $Y\setminus D_{Y}$, and the inverse image of each cusp of $Y$ is the elliptic tail attached at the corresponding point of $D_Z$. If $Y$ and $Y'$ are pseudostable curves whose $Z$ are the same, then $Y$ and $Y'$ are isomorphic. 
\end{corollary}

Next, Schubert considers flat families over $B=\spec(R)$ with smooth connected general fibers of genus $g \ge 3$ and reduced special fibers as shown in the diagrams below:
\newdiagramgrid{dvrgrid}{0.5,0.6,0.5,0.6, 0.5, 1.4,0.5,0.6,0.5,0.6, 0.5,}{0.5,0.6, 0.5,0.6,0.5,}
\setdiagram{schubertdvrdiagram}{nohug,grid=dvrgrid}{
 	&	&Z	&	&	&	&	&	&Z_{\eta}	&	&	\\
	&	\ldTo^{\pi_{X}}	&	&\rdTo^{\pi_{Y}}	&	&	&	&\ldTo^{\pi_{X,\eta}}	&	&\rdTo^{\pi_{Y,\eta}}	&	\\
X	&	&	&	&Y	&	&	X_{\eta}&	&\rTo^{\phi}_{\cong}	&	&Y_{\eta}	\\
	&	\rdTo_{\psi_X}&	&\ldTo_{\psi_Y}	&	&	&	&	\rdTo_{\psi_{X,\eta}}&	&\ldTo_{\psi_{Y,\eta}}	&	\\
	&	&B	&	&	&	&	&	&B_{\eta}	&	&	\\
	}%

In the situation of~\eqref{schubertdvrdiagram}, call $p\in X_0$ exceptional if $\pi_{X,0}$ is not an isomorphism over any open neighborhood of $p \in X_0$. When an exceptional $p$ is fixed, define $E=\pi_{X,0}^{-1}(p)$, $F=\overline{Z_0\setminus E}$, $C=\pi_{Y,0}(E)$, $D= \overline{Y_0\setminus C}$ and $l=\#(C\cap D)$.

\beginS{lemma} \label{schubertclaims}
Given $X\to B$, $Y\to B$ and $\phi$ as in~\eqref{schubertdvrdiagram}, there is a flat $Z$ satisfying the following:
\beginS{enumerate}
\item \label{ciso} $\pi_{X,\eta}$ is an isomorphism.
\item \label{cred} $Z_0$ is reduced and no component of $Z_0$ is collapsed by both $\pi_{X,0}$ and $\pi_{Y,0}$.
\item \label{ceandf} If $p$ is exceptional, then $\pi_{Y,0}$ is an isomorphism over every point of $C$ except those in the image of $E\cap F$.
\item \label{czero} If $p$ is an exceptional smooth point or node of $X_0$, then $E$ is a curve of genus $0$.
\end{enumerate}
\end{lemma}
	
\newcommand{\cref}[1]{(\ref{schubertclaims}.\ref{#1})}	
For these claims, in which pseudostability plays no part, I refer to~\cite{Schubert}*{Lemmas 4.4 to 4.7}.

The application to moduli then follows in his Lemma~4.8 that asserts,
\beginS{lemma}\label{dmversuspseudostable}
If, in the situation above, $X_0$ is Deligne-Mumford stable and $Y_0$ is pseudostable, then:
\beginS{enumerate}
	\item $\pi_{X}:Z\to X$ is an isomorphism.
	\item $\pi_{Y}:Z\to Y$ is an isomorphism except over cusps of $Y_0$.
	\item The inverse image under $\pi_Y$ of a cusp in $Y_0$ is an elliptic tail in $Z_0$.
\end{enumerate}
\end{lemma}
\beginS{proof}
To prove the first claim, it suffices to show that there are no exceptional points $p$. If $p$ is exceptional, then by~\cref{czero} $E$ is a curve of genus $0$ meeting $F$ in either $1$ or $2$ points (since $p$ is at worst nodal) and every point of $C\cap D$ is the $\pi_{Y,0}$-image of one of these by~\cref{ceandf}, so $l \le 2$.

If $l=0$, then $Y_0 =C$ and, by~\cref{ceandf}, there are $n\le 2$ points at which $\pi_{Y,0}:E\to C$ is not an isomorphism. Checking that $C$ has genus $n$ for each possible value, we get a contradiction to $g\ge 3$.

If $l=1$, then $C$ is smooth at this point and $\pi_{Y,0}:E\to C$ is an isomorphism near it. There is at most one other point in $E\cap F$ and hence, by~\cref{ceandf}, at most one point near which $\pi_{Y,0}:E\to C$ is not an isomorphism. If there are none, then $C$ is a rational subcurve of $Y_0$ meeting the rest of $Y_0$ in at most $2$ points. If there is $1$, this point must be a cusp and $C$ must be a rational cuspidal tail. In both cases this contradicts pseudostability of $Y_0$.

If $l=2$, then $C$ is again non-singular at these points and hence 
$\pi_{Y,0}:E\to C$ is an isomorphism near them and hence, by~\cref{ceandf}, everywhere. Once again $C$ is a rational subcurve of $Y_0$ meeting the rest of $Y_0$ in at most $2$ points.

Next note that the notion of exceptional point makes equal sense for points $q$ of $Y_0$ and that claim (4) of Lemma~\ref{schubertclaims} holds equally for these. Applied to an exceptional point $q$ of $Y_0$ that is not a cusp, it shows that the inverse image $E$ of $q$ in $Z_0$ will be a genus $0$ subcurve meeting its complement in $Z_0$ at $1$ point if $q$ is smooth and in at most $2$ points if $q$ is a node. This contradicts the first claim and proves the second.

Finally, if $q\in Y_0$ is a cusp, then it must be exceptional, as $Z_0$ has no cusps, and must have a unique preimage on $F$ that is a smooth point of $F$. Thus $F$ must be the normalization of $Y_0$ at $q$ and have genus $(g-1)$. But then, since $E\cap F$ is a single point, $E$ must be a connected curve of genus $1$ and the last claim follows.
\end{proof}

The remainder of Schubert's argument closely follows the outline in subsection~\ref{pluricanonicalproperties}. The one point of difference is in the proof of the analog to Corollary~\ref{StabilityAndJss}.\ref{semistablereductionstep} showing that $3$- and $4$-canonical models of pseudostable curves are asymptotically Hilbert stable. It is here that the Lemmas above are invoked. Fix $\nu$ to be 3 or 4.

Given a pseudostable curve $Y'_0$, let $Z_0$ be the stable curve with standard tails of Corollary~\ref{stablemodelofpseudo} and let $Z$ be a flat smoothing of $Z_0$ over $B$ and $ \alpha:B\to\h $ be the map induced by taking the family of $\nu$-canonical models of $Z$. Since the generic fiber $ Z_{\eta} $ is smooth, its Hilbert point $ [Z_{\eta}] $ lies in $\semi{\h}$ by Theorem~\ref{StabilityOfSmoothCurves}. 

After possibly making a finite change of base, we can find a map $ \beta:B\to\semi{\h} $ that agrees with $ \alpha $ at $ \eta $. By pulling back the universal curve over $\h$\ by $ \beta $, we obtain a second curve $ Y\to B $ whose generic fiber is also $ Z_{\eta} $. Since $Y_0$ is Hilbert semistable, Schubert's Potential Stability results imply that $Y_0$ must be a pseudostable curve. But we are also in the situation of Lemma~\ref{dmversuspseudostable} and last claim there implies that the $Z_0$ is also the stable curve with standard tails determined by $Y_0$. Corollary~\ref{stablemodelofpseudo} then implies that $Y'_0\cong Y_0$ so is semistable. Strict semistability is ruled out as in the proof of (\ref{StabilityAndJss}.\ref{semistablereductionstep}) by the fact that any Hilbert semistable curve has finite automorphism group.


\section{Weighted Pointed Curves}
\label{swinarski}

In this section, I want to review the recent \GIT construction by Swinarski of the moduli spaces $\mgnbar$ of $n$-pointed stable curves of genus $g$. With minor adjustments, the proof also constructs the moduli spaces of weighted pointed curves of Hassett~\cite{HassettWeighted}. Here I will first explain what's involved in setting up the GIT problem and give analogues of the criteria of Section~\ref{numerical}. Then I'll sketch the ideas in the proof that smooth pointed curves have stable Hilbert points with respect to suitable linearizations. I will explain why Gieseker's Criterion is not an adequate tool here and prove the main result, the Span Lemma~\ref{spanlemma} that substitutes for and sharpens it. However, I will omit the essentially combinatorial verification that it suffices to check stability, giving only a statement of the main result.  Discussion of the modifications need to get a Potential Stability theorem and construct moduli for stable pointed curves  is postponed those until the next section on stable maps, since the similar issues arise in both constructions.

\subsection{Setting up the \GIT Problem}
\label{swinarksisetup}

To begin, we need to set up parameter spaces for pointed curves, describe suitable linearized line bundles on them and understand the numerical criterion in terms like those in Section~\ref{numerical}. We describe how to do this in the notation established there. 

Let $X = (C, [p_1, \ldots, p_n])$ be an curve of genus $g$ and degree $d$ in $\PV$ and an ordered set of $n$ points, not yet necessarily distinct or lying on $C$. Viewing $\PV$ as the Hilbert scheme of points, i.e. subschemes with Hilbert polynomial (and function) the constant $1$, $X$ determines a point in the product $\h^* := \h \times \prod_{k=1}^n \PV$. For $m$ sufficiently large and $m_k>0$ for each $i$, we fix $L:= L(m; m_1, \ldots,m_n)$ to be the very ample line bundle obtained by tensoring the pullback $\Lambda_m$ from $\h$---taken as in Section~\ref{numerical} with respect to its Pl\"ucker embedding in $\PP(W_m)$---with the pullbacks of $\o_{\PV}(m_k)$ from the $\thst{k}{th}$ factor $\PV$, for all $k$. Since each of the factors is naturally linearized, so is $L$. We will call this the linearization with parameters $(m; m_1, \ldots,m_n)$ and will fix it henceforth. 

We will, once again, express the numerical criterion in terms  this linearization of $L$. We can simplify notation slightly by dropping the trivial first exterior power from the point factors of the space $W':= W(m; m_1, \ldots,m_n)$ on which $\ssl(V)$ naturally acts and writing 
$$W' = W_m \otimes S_{m_1} \otimes \cdots\otimes S_{m_n}$$
If we then fix $1$-ps $\rho$ as in~\eqref{setuponeps}, coordinates $Z$ that diagonalize the action of $\rho$ on $W$ are determined by the data of a set $z$ of $\p(m)$ degree $m$ monomials in the $x_i$, and, for each $k$, a single monomial $y_k$ of degree $m_k$. A coordinate $Z$ is non-zero at the point $[X] \in \h^*$ if and only if the elements in $z$ restrict on $C$ to a basis of $\Hzeroxom{C}$, as before, and each $y_k$ is non-zero at $p_k$ (the analogous condition for the point factors). 

Define $w_B(k)$ to be the least weight of a coordinate $x_i$ not vanishing at $p_k$. Clearly, the least weight of a $y_k$ non-zero at $p_k$ is then $m_k w_B(k)$. Thus the least weight of a coordinate $Z$ non-zero at $[X]$ is $w_B(m)+ \sum_{k-1}^n m_k w_B(k)$. This gives the first claim below from which the second is immediate, setting $w_F(k) = w_B(k)$.

\beginSS{proposition}{[Numerical criterion for pointed curves]} \label{NumericalCriterionForPointedCurves} 
A point $[X]\in\h^*$ is $(m; m_1, \ldots,m_n)$-stable \resp{semistable} if and only if the equivalent conditions below hold:
\beginS{enumerate}
	\item For every weighted basis $B$ of $V$, there is a $B$-monomial basis of \\ $\Hzeroxom{C}$ such that $w_B(m)+ \sum_{k-1}^n m_k w_B(k) < \text{\resp{$\le$}~} 0$.
	\item \label{filtrationversion} For every weighted filtration $F$ of $V$ whose weights $ w_{i}$ have average~$\alpha$, $$ w_{F}(m)+ \sum_{k-1}^n m_k w_F(k) < \text{\resp{$\le$}~} (m\p(m) + \sum_{k=1}^n m_k)\alpha\,. $$
	\end{enumerate}
\end{proposition}

To ensure that the points $p_k$ lie on the curve $C$, we simply need to replace $\h^*$ by the closed subscheme $\hhat$ determined by the  incidence conditions $p_k \in C$, for all $k$:  $\hhat$ is defined scheme theoretically by the condition that the ideal of $C$ is in the ideal of each $p_k$. 

To get interesting consequences, we need to balance $m$ and the $m_k$'s. As we scale the former to get asymptotic results, we want the latter to scale correspondingly. 
\beginS{definition}\label{blinearization}
For $B=(b_1,\ldots, b_n)$ with each $b_k>0$, the $B$-linearization of the $\ssl(V)$ action on $\h^*$ is that given by setting $m_k = b_k \frac{m^2}{2}$ and it is then convenient to set $b:= b_B=\sum_{k=1}^n b_k$.  We can, and will, allow $b_k$ to be rational, understanding that we always take $m$ sufficiently divisible that all the $m_k$ are integral. 
\end{definition}

I should note that my $b$'s are twice those in Swinarski~\cite{SwinarskiThesis}---I'll explain this change in a moment. However, the $m^2$ factor in the definition of $b_k$ is fundamental. It is needed to make the ``point'' terms in Proposition~\ref{NumericalCriterionForPointedCurves} have the same order in $m$ as the ``curve'' terms and makes possible an immediate analogue of Lemma~\ref{AsymptoticNumericalCriterion}.

\beginSS{lemma}{[Asymptotic numerical criterion for pointed curves]} \label{AsymptoticNumericalCriterionForPointedCurves} 
Let $X$ be an $n$-pointed curve of degree $d$ and genus $g$ in $\PV$. For a weighted filtration $F$, let $b_F := \sum_{k=1}^n b_k w_F(k)$.
	\beginS{enumerate} 
		\item \label{ealphaInequalitypointedcurves} If $e_{F} + b_F < 2\alpha_{F}(d+b)$, then $X$ is $B$-Hilbert stable with respect to $F$. 			
		\item \label{ealphaInequalitypointedcurvesunstable} If $e_{F}+ b_F> 2\alpha_{F}(d+b)$, then $ X$ is $B$-Hilbert unstable with respect to~$F$. 
		\item \label{asymptoticuniformpointedcurves} If there is a $\delta>0$ such that $$ e_{F}+b_F < 2\alpha_{F}(d+b)-\delta $$ for all weighted filtrations $F$ associated to the Hilbert point of any $X$ in a subscheme $\sss$ of $\hhat$, then there is an $M$, depending only on $\sss$, such that the \thst{m}{th} Hilbert point $[X]_{m}$ of $X$ is $B$-stable for all $m\geq M$ and all $X$ in $\sss$. 
	\end{enumerate} 
\end{lemma}

We think of $b_k$ as a weight on the point $p_k$ in the sense of Hassett~\cite{HassettWeighted} and, with the chosen scaling, stability forces $b_k \in (0,1]$ and more. To see this, and get a feel for stability for pointed curves, let's work an example, assuming that $g\ge 2$ to simplify. Consider the filtration $F$ that assigns weight $0$ to sections vanishing at smooth point $q$ of $C$ and weight $1$ to all others. First, for $i=0, \ldots, m$ there is a section in $\Hzeroxom{C}$ vanishing to order exactly $i$ at $q$ and having weight $m-i$ so, summing over $i$, $w_{F}(m) = \frac{1}r{2}m^2 + O(m)$ and $e_F=1$. Each $w(k)$ is either $1$ or $0$ according as $p_k$ is or is not equal to $q$ so $b_F=\sum_{p_k=q}b_k$.

The right hand side of~\ref{AsymptoticNumericalCriterionForPointedCurves}.\ref{ealphaInequalitypointedcurves} is $\frac{d+b}{N}$, which, if we let $d\to \infty$ approaches $1$ from above. Thus if $b_F \le 1$, $X$ is $F$-stable and if $b_F >1$, then, for large $d$, $X$ is $F$-unstable. Thus, taking $q = p_k$, the weight of any marked point $p_k$ on a stable $X$ can be at most $1$ when $d \gg 0$, and if so, no other marked point can equal $p_k$. More generally, semistability implies that the sum of $b_k$ over the set of marked points that are equal to $q$ can be at most $1$, exactly Hassett's condition. This is my justification for choosing the normalization I do of the $b_k$s: Swinarski's can be at most $\frac{1}{2}$.

I'll leave the reader to check similarly that if $q$ is a node of $C$, then $e_F=2$ and hence, when $d \gg 0$, a stable $X$ must have $b_F=0$. In other words, no marked point can be a node of $C$.

\subsection{Hilbert stability of smooth pointed curves}
\label{swinarksismooth}

In this subsection, I want to explain the new geometric ideas behind the estimates that Swinarski uses to prove Hilbert stability of smooth pointed curves, for which Gieseker's Criterion~\ref{GiesekersCriterion} turns out to be insufficiently sharp and that are the main novelties in his construction. Even with these in hand, the combinatorial  argument deducing stability is delicate and lengthy (over 20 pages) so I'll simply outline the strategy and give some motivating examples, referring to ~\cite{SwinarskiThesis} for the details.  

To get a feel for the difficulties, let's first look at his Example 1 (defined in 2.4 and discussed at several subsequent points). We let $n=3$ $L=\o_C(1)$ and consider the filtration $F$ for which there are four weight spaces $V_0\supsetneqq V_1\supsetneqq V_2\supsetneqq V _3$ given by, 
{\stretcharray{1.2}{0pt}{
\beginS{equation}
\label{exampleone}
\begin{array}{c@{~}l@{~}c@{~}l@{~}c@{~}l@{~}c}
 V&\supsetneqq&\HO{C}{L(\text{--}p_1)}&\supsetneqq&\HO{C}{L(\text{--}p_1\text{--}p_2)}&\supsetneqq&\HO{C}{L(\text{--}p_1\text{--}p_2\text{--}p_3)}\\
w_{0}=\frac{3}{6}&>&w_{1}=\frac{2}{6}&>&w_{2}=\frac{1}{6}&>&w_{3}=0 
\end{array}
\end{equation}
}}%
This has normalized weights (decreasing to $0$ and summing to $1$), hence $\alpha_F=\frac{1}{N}$ and, in Lemma~\ref{AsymptoticNumericalCriterionForPointedCurves}, the right hand sides are slightly larger than $2$ for large $d$. Since $e_i=i$ for $i \le 3$ and $e_i=3$ for $i >3$, it's easy to  see that the $\epsilon_F$ of Gieseker's Criterion~\ref{GiesekersCriterion} equals $\frac{3}{2}$ (achieved by any subsequence). The least weights of sections not vanishing at $p_k$ are $\frac{3}{6}, \frac{2}{6}, \frac{1}{6}$ for $k=1,2,3$ so $b_F=1$. Thus if we estimate $e_F$ by $\epsilon_F$, then we only get $e_F + b_F \le \frac{5}{2}$. Below, we'll see that Swinarski's estimates show that $e_F$ is really $1$. Hence $e_F + b_F =2$ and $X$ is $F-$stable. But we can already see that substantially better estimates than those used in Gieseker's criterion are needed to deduce stability for pointed curves.

The key idea can be understood by looking a bit further at this example. Recall that the basic idea in Gieseker's criterion is to combine estimates of the codimensions of various weight spaces in $\HO{C}{\lbpow{L}{m}}$. One of these estimates uses $\HO{C}{\lbpow{L}{m}(-mp_1)} = V_1^m$ to estimate by $m$ the codimension of the space of sections of weight at most $\frac{1}{3}m$. But $\HO{C}{\lbpow{L}{m}(-\frac{1}{2}mp_1-\frac{1}{2}mp_2 )} = \bigl(V_0 V_2\bigr)^{\frac{1}{2}m}$ also has codimension $m$ and weight at most $\frac{1}{3}m$. These two spaces intersect in the space
$\HO{C}{\lbpow{L}{m}(-mp_1-\frac{1}{2}mp_2 )}$ of codimension $\frac{3}{2}m$ so their span is a subspace of weight at most $\frac{1}{3}m$ with codimension $\frac{1}{2}m$. In other words, the codimension estimate for the $\frac{1}{3}m$-weight space used in Gieseker's criterion in this case is a factor of $2$ from being sharp. Plugging this improved estimate into~\eqref{singleFestimate} reduces $e_F$ by $\frac{1}{3}$.

Swinarski's strategy is to use this idea to improve codimension estimates in a systematic way and see that the gains made suffice to prove Hilbert stability. Essentially, he replaces ``monomial'' subspaces (subspaces like $V_1^m$ that are spanned by monomials in the $x_i$) by ``polynomial'' subspaces (spans of several mononial subspaces of the same maximum weight like $V_1^m$ and $\bigl(V_0 V_2\bigr)^{\frac{1}{2}m}$). The key estimate is:

\beginSS{lemma}{[Span Lemma]} \label{spanlemma} Fix a smooth curve $C$ of genus $g$ embedded in $\PV$ by a line bundle $L$ of degree $d$. Fix a set $\{q_1, \ldots, q_a\}$ of points of $C$, and an $A \times a$ matrix $[c_{i,j}]$ of non-negative integers. For each $j$, let $\cmin_j$ and $\cmax_j$ be the minimum and maximum entries in the $\thst{j}{th}$ column. For $i=1, \ldots, A$, define $E_i := \HO{C}{\lbpow{L}{m}(-\sum_{j=1}^a c_{i,j}q_j)}$ and define $E:= \mcspan(E_1, \ldots, E_A)$. 

If $\sum_{j=1}^a \cmax_j <dm-2g$, then 
$E =  \HO{C}{\lbpow{L}{m}(-\sum_{j=1}^a \cmin_j q_j)}$ 
and $E$ has codimension $\sum_{j=1}^a \cmin_j$ in $\HO{C}{\lbpow{L}{m}}$.
\end{lemma}
The example in~\eqref{exampleone} is a toy case of the Span Lemma with $A=a=2$, $q_1=p_2$ and $q_2=p_1$. In more general applications, the $q_j$'s may or may not be marked points on $C$. The proof will show that any $E_i$ for which no $d_{i,j}=\cmin_j$ can be discarded without affecting $E$ so in applications we can always take $A \le a$.
\beginS{proof}
We first verify the codimension of $E$. If $I$ is a non-empty subset of $\{1, 2, \ldots, A\}$, let $E_{I}= \cap_{i\in I}E_i$.  For any subspaces $E_i$ of a vector space $W$, 
\beginS{equation*}
\label{subspaceinex}
\codim\bigl(\mcspan(E_i, \ldots, E_A)\bigr) = \sum_{I} (-1)^{(\#I-1)}\codim(E_I)
\end{equation*}
On the other hand, $E_{I}= \HO{C}{\lbpow{L}{m}(-\sum_{j=1}^a c_{I,j}p_j)}$ where $c_{I,j}= \max_{i \in I} c_{i,j}\le \cmax_j$ and the hypothesis on the $\cmax_j$ lets us compute the codimension of $E_{I}$ by Riemann-Roch as $\sum_{j=1}^a \cmax_{I,j}$. Plugging this into~\eqref{subspaceinex} gives
\beginS{equation*}
\codim\bigl(\mcspan(E_i, \ldots, E_A)\bigr) 
= \sum_{I} (-1)^{(\#I-1)} \sum_{j=1}^a c_{I,j} = \sum_{j=1}^a \sum_{I} (-1)^{(\#I-1)}  c_{I,j} \,. 
\end{equation*}
Fix $j$ and assume the $c_{i,j}$ are distinct. Then every $c_{I,j}$ is a $c_{i,j}$ for a unique $i \in I$ and the number of $I$ for which $c_{i,j}$ gives $c_{I,j}$ is $2^k$ where $k= \#\{i'| c_{i',j} < c_{i,j}\}$. Moreover, except when $k=0$, exactly half of these $I$ are of each parity and so cancel in the $\thst{j}{th}$ summand above. Making these substitutions and cancellations, the only term remaining in the $\thst{j}{th}$ summand is that where $I = \{i\}$ and $c_{i,j} = \cmin_j$. The same argument applies when the $c_{i,j}$ are not distinct, if we ``break ties'' by replacing $<$ with any total order refining $\le$.

To see that $E =  \HO{C}{\lbpow{L}{m}(-\sum_{j=1}^a \ct_jq_j)}$, it now suffices, since the codimensions are equal, to show either containment. But every section in each $E_i$ vanishes to order at least $\cmin_j$ at $q_j$ so the same is true of every section in their span $E$.
\end{proof}

The combinatorial argument by which Swinarski deduces stability does not involve a reduction to a linear programming problem. Instead, he works with the fixed subsequence of the $V_j$ at which the base locus jumps and estimates separately the codimensions of each of the stages in Gieseker's double filtration~\eqref{giesekerfiltration} using the Span Lemma~\ref{spanlemma}. The basic idea is as follows. Fix one of Gieseker's subseries $U_{k,l}$ and a point $q_i$ in the base locus of $V_N$. Then, find a subseries $U_{k,l,i}$ of $U_{k,l}$ obtained by restricting a product of the form $\sym^{n}\bigl(V\cdot\sym^{(p-w)}(V_{s}) \cdot\sym^{w}(V_{t})\bigr)$ 
with $s$, $t$ and $w$ chosen so that, first, the maximum weight of a section in $U_{k,l,i}$ is no greater than for $U_{k,l}$ and, second, the multiplicity of $q_i$ in the base locus of $U_{k,l,i}$ is as small as possible consistent with this weight requirement. (We can always take $U_{k,l,i} = U_{k,l}$ if necessary.) Then let $U'_{k,l}$ be the span of all the $U_{k,l,i}$ and use the double filtration given by these to estimate $w_F(m)$. The Span Lemma quantifies the drop in codimension from each $U$ to the corresponding $U'$.

Swinarski works out these estimates in two stages. First, he computes a \emph{virtual profile} that graphs the points $\bigl(\codim(U'_{k,l}), \textrm{weight}(U'_{k,l})\bigr)$ that would arise if \emph{fractional} $w$ could be used above. This is a piecewise linear function with one segment for each $k$ in the subsequence. The virtual profile that arises from the example defined in~\eqref{exampleone} is shown in Figure~\ref{ExampleOneDiagram}. The arrows show, on each segment, the virtual base locus of the $U'$ with virtual weight $\frac{1}{2}m-\alpha$. For example, the top segment, of slope $-2$, is what results from the discussion following~\eqref{exampleone}. Similar considerations involving the other $p_k$ lead, via the Span Lemma, to the other segments. Note that the area under this profile is $\frac{1}{2}$, just what is needed to check stability in this example. 

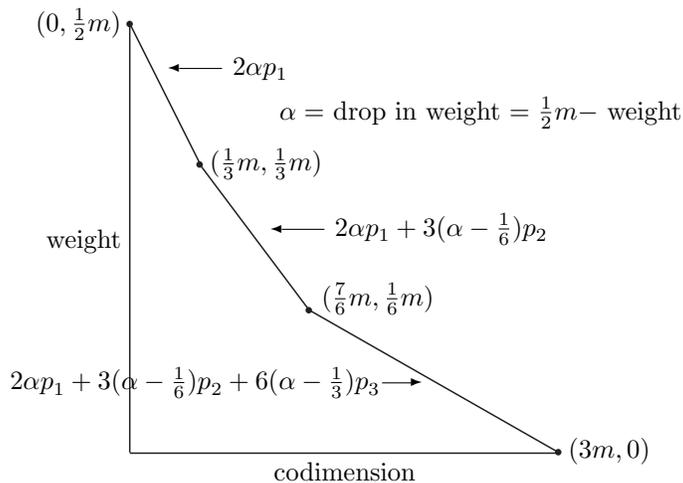
\capdraw{ExampleOneDiagram}{Virtual profile for the example of \eqref{exampleone}.}{.75}{ExampleOne.pdf}{
\etext{.02}{.97}{$(0,\frac{1}{2}m)$}
\wtext{.21}{.67}{$(\frac{1}{3}m,\frac{1}{3}m)$}
\wtext{.45 }{.38}{$(\frac{7}{6}m,\frac{1}{6}m)$}
\wtext{.98}{.05}{$(3m,0)$}
\wtext{.10}{.87}{$\lTo 2\alpha p_1$}
\wtext{.32}{.53}{$\lTo 2\alpha p_1+ 3(\alpha-\frac{1}{6})p_2$}
\etext{.66}{.2}{$2\alpha p_1+ 3(\alpha-\frac{1}{6})p_2+ 6(\alpha-\frac{1}{3})p_3\to$}
\wtext{.36}{.78}{$\alpha = $ drop in weight $=\frac{1}{2}m -$ weight} 
\ntext{.5}{.02}{codimension}
\etext{.02}{.5}{weight}
}%

Unfortunately, the area under virtual profile does not bound $w_F(m)$ because of the need to choose integral $w$. Using these $w$ produces a \emph{profile} that is the graph of a step function with $m$ steps over each segment in the virtual profile and that lies above both the virtual profile and the graph of the step function determined by the weight filtration. The area under this profile \emph{does} bound $w_F(m)$.

The combinatorial task is then to estimate both the area under the virtual profile and the area between the profile and the virtual profile and to show the resulting estimate for $e_F$ is sharp enough to yield Hilbert stability. Defining the profile and virtual profile in general and carrying this out occupies pages 15-38 of the online version of~\cite{SwinarskiThesis}. The main result is:

\beginS{theorem}\label{stabilityofpointedsmoothcurves} 
Fix $d$ and $g$ and let $\h^*$ be the corresponding Hilbert scheme of $n$-pointed curve of genus $g$ in $\PP^{N-1}$. Fix a linearization $B$, as in~\eqref{blinearization}, satisfying $b_k <1$ for all $k$ and, for some positive $\delta$,  $b = \sum_{i=1}^n b_k = \frac{g-1}{N-1} + \delta$. Let  $X=(C, [p_1, \ldots, p_n])$ be any smooth $n$-pointed curve such that for any point $q$ of $X$, $\sum_{p_k=q} b_k \le 1$. Then $X$ is asymptotically $B$-Hilbert stable. 
\end{theorem}

A few remarks are in order. This is the main Case A of Swinarski's Theorem~8.2.1. He handles some other linearizations in other cases. Recall also that his $b_k$ are half mine (and my $\delta$ is his $\epsilon N$). Finally, he only claims $m$-Hilbert stability for a certain infinite set of $m$. However, his argument shows that the polynomial representing $w_F(m)$ has negative leading coefficient bounded away from $0$ in terms of $\delta$ and that $w_F(m)$ has negative value for an $m$ only depending on the global parameters $d$, $g$, $B$ and $\delta$. From this, we may deduce a uniform upper bound for the linear coefficient in $w_F(m)$ and hence its negativity for all sufficiently large $m$.


\section{Pointed Stable Maps}
\label{baldwin}

In this section, I want to review the recent \GIT construction by Baldwin and Swinarski~\cite{BaldwinSwinarski} of the moduli spaces $\mgnbardhat$ of $n$-pointed stable maps of genus $g$ with image of degree $\dhat$ in $\PP^{r}$ for all $(g,n)$ except $(0,0)$, $(0,1)$ and $(1,0)$. This construction is not self-contained but depends on the existence of $\mgnbardhat$ as a coarse moduli space due to Fulton and Pandharipande~\cite{FultonPandharipande}. However, Baldwin~\cite{BaldwinPositive} has modified their argument to give a standalone proof that is even valid in those positive characteristics not dividing $\dhat$ for which the corresponding moduli problem is separated.

Once again, the proof is far too long to give in full. I have not reproduced even the statements of many of their results because the quantifications alone often stretch to a dozen lines. Instead, I have chosen to sketch the main points at which it diverges from the model in Section~\ref{gieseker}, sending the interested reader to~\cite{BaldwinSwinarski} for almost all details. The main novelty is Baldwin's ingenious inductive argument for the stability of pointed maps with smooth sources. This is treated in the second subsection. The first deals with the setup of the \GIT problem they treat and the new issues that arise in proving a Potential Stability theorem for pointed maps. In particular, it explains why an appeal to another construction is needed. 

\subsection{Overview of the construction for pointed stable maps}
\label{overviewformaps}

Setting up parameter spaces, fixing suitable linearized line bundles on them and interpreting the numerical criterion are all straightforward variants of Section~\ref{swinarski}. We describe briefly how to do this in the notation established there. 

Set $Z = \PV \times \PP^{r}$ with $\pi$ and $\hat\pi$ the two projections and let $\o_Z(m,\mhat)$ be the restriction to $Z$ of $ \pi^*\bigl(\o_{\PV}(m)\bigr) \otimes  {\hat\pi}^*\bigl(\o_{\PP^{r}}(\mhat)\bigr) $.  Then let $\h'$ be the Hilbert scheme of subcurves $C'$ of $Z$ of bidegree $(d,\dhat)$, that is, having Hilbert polynomial $P(m,\mhat)= dm+\dhat \mhat+(1-g)$. The natural models of $\h'$ are constructed by restriction of sections of $\o_Z(m,\mhat)$ for suitably large $m$ and $\mhat$ and have \plucker coordinates indexed by sets of $P(m,\mhat)$ polynomials of bidegree $(m,\mhat)$ corresponding to bundles $\Lambda_{m,\mhat}$. Note that while the linear series cut out by forms of bidegree $(1,1)$ on $C'$ is tautologically very ample, it does not follow that either of the linear series induced on $C'$ by projecting $Z$ to $\PV$ or $\PP^{r}$ are. 

Let $X = (C, [p_1, \ldots, p_n])$ be a prestable $n$-pointed curve of genus $g$ (meaning that $C$ reduced, connected, complete and with only nodes and that the marked points are distinct, smooth points of $C$). Let $f:C\to \PP^{r}$ be a stable map with image of degree $\dhat$. Stability for $f$ means finiteness of automorphisms of $X$ commuting with it, which rules out $(g,n)= (1,0)$ and in all other cases means that any smooth rational curve collapsed by $f$ must contain at least $3$ special points. Fix also an embedding of $C$ in $\PV$ as a curve of degree $d$. Then the graph $\Gamma$ of $f$ in $ \subset \PV \times \PP^{r}$ determines a point $C'$ of $\h'$. In this case, the linear series induced on $C'$ by projecting to $\PV$ is necessarily very ample but that given by projection onto $\PP^{r}$ need not be.

We can record the $n$-marked points of $C$, as before by taking a product of $\h'$ with $n$ copies of $Z$ and passing to the incidence correspondence $\hhat$ for which these points lie on $\Gamma$. We fix $L:= L(m, \mhat; m_1, \ldots,m_n)$ to be the very ample line bundle obtained by tensoring the pullbacks of the corresponding $\Lambda$'s from each factor. In~\cite{BaldwinSwinarski}, the $m_k$ are always taken to have a common value $m'$ so we will write $L:= L(m,\mhat, m')$. I have also suppressed the $\mhat_k$'s that might be formally expected in the point factors because they will play no role in the sequel.

The reason for this is that there is only one way for $\ssl(V)$ to act on $\PP^r$---trivially, so that $1$-ps's (or filtrations) always act with $0$ weights on the $\PP^r$ side of the picture. We then extend all the notions defined in terms of these weights in the obvious bigraded way. For the point factors $w_B(k)$, we can thus drop the dependence on $\mhat$ since we can always find some monomial on the right side non-vanishing at $p_k$. 

\beginSS{proposition}{[Numerical criterion for stable maps]} \label{NumericalCriterionForStableMaps} 
A point $[f] \in\hhat$ is stable \resp{semistable} with respect to the linearization defined above with parameters $(m, \mhat, m')$ if and only if the equivalent conditions below hold:
\beginS{enumerate}
	\item For every weighted basis $B$ of $V$, there is a $B$-monomial basis of \\ $\HO{C}{\o_C(m, \mhat)}$ such that 
	$$w_B(m, \mhat)+ m' \sum_{k=1}^n w_B(k) < \text{\resp{$\le$}~} 0\,.$$
	\item \label{filtrationversionstablemaps} For every weighted filtration $F$ of $V$ whose weights $ w_{i}$ have average~$\alpha$, $$ w_{F}(m, \mhat)+ m' \sum_{k=1}^n  w_F(k) < \text{\resp{$\le$}~} \bigl(m\p(m,\mhat) + n m'\bigr)\alpha\,. $$
	\end{enumerate}
\end{proposition}

Since Baldwin and Swinarski work mainly with the criterion above, 
I leave the reader to formulate the obvious asymptotic variant.

Next they define $J := J_{\nu,c}$ to be the locus in $\hhat$ where $C$ is prestable, the linear series given by $\PV$ is non-degenerate and  the associated invertible sheaf $\o_C(1, 0)$ is isomorphic to $\lbpow{\bigl(\omega_C(\sum_{k=1}^n p_k)\bigr)}{\nu}\otimes \o_C(0, c\nu)$. (They use $a$ for my $\nu$.) This $J$ is locally closed in $\hhat$ (cf.~\cite{FultonPandharipande}*{Proposition 1}). By the universal property of $\hhat$, $J$ has the local universal property for the moduli problem of $n$-pointed stable maps of genus $g$ with image of degree $\dhat$ in $\PP^{r}$ and, by construction, $\ssl(V)$ orbits in $J$ correspond to isomorphism classes of such maps. To see that $J \quotient \ssl(V)$ is a coarse moduli space for this problem, it then suffices to find an very ample linearized $L$ such that $\semi{J}(L) = \stable{J}(L) = J$ where $\semi{J}(L)$ and $\stable{J}(L)$ be the subsets of $J$ of points  semi-stable and stable, respectively, with respect to the chosen linearization $L$ (cf.~\cite{NewsteadModuli}*{Proposition~2.13}).  In fact, it suffices to find an $L$ for which $\semi{J}(L) = J$ since all points of $J$ have finite stabilizers. 

To show the first containment, $\semi{J}(L) \subset J$ for suitable $L$, Baldwin and Swinarski prove a Potential Stability Theorem (their Theorem~5.19). Its hypotheses involve a complicated set of inequalities on the invariants $g, \nu, c$ and $\dhat$ of the basic setup, on quantities like the pluricanonical $d$ and $N$ determined by these, on further quantities defined in terms of these in Proposition~4.6 as uniform bounds on $\hhat$ of various sorts, and, finally, on the parameters $(m,\mhat, m')$ of the linearization used. 

Rather than give these in detail, I will simply note that they determine a open convex set $\mathbf{L}$ of linearizations and hold whenever $\nu \ge 10$ (and in many cases for smaller $\nu$), as long as $m$ is sufficiently large and $|\frac{\mhat}{m} - \frac{c\nu}{2\nu-1}|$ and $|\frac{m'}{m^2} - \frac{\nu}{2\nu-1}|$ are both sufficiently small (cf. the Remark following Theorem~5.21). That is, even more than for weighted pointed curves, it is important to carefully balance the curve, map and point parameters $m$, $\mhat$ and $m'$ of the linearization. In particular, for large $m$, $\mathbf{L}$ always contains the \emph{central linearization} $L^{*}$ for which $\mhat=\frac{c\nu}{2\nu-1}m$ and $m'= \frac{\nu}{2\nu-1}m^2$.

Theorem 5.19 adds, naturally, the conditions that the marked points be distinct and smooth on $C$ to the notion of potentially stable in Definition~\ref{PotentiallyStable}. A more important difference is that it uses a more complicated version of the Subcurve Inequality~\ref{PotentiallyStable}.\ref{SubcurveInequality}---(iv)  in  Theorem~5.19. 

There are some significant new difficulties to overcome in the proof. The most substantial come at the start, when it is necessary to consider $C \subset Z$ for which it is not known that projection to the $\PV$ factor is of degree $1$. To establish this, Baldwin and Swinarski first prove that $\pi_{\PV}(C)$ is nondegenerately embedded in $\PV$.  Next, they show that $\pi_{\PV}(C)$ is generically reduced. Both these arguments follow the lines of Gieseker's proofs for the same assertions. In the context of maps, the second argument shows, at the same time, that $\pi_{\PV}$ must have degree $1$ over every component of $\pi_{\PV}(C)$: it cannot send two components of $C$ to the same component $\pi_{\PV}(C)$ nor map any component of $C$ multiply to its image in $\pi_{\PV}(C)$.  

It thus remains to show that no component of $C$ collapses to a point of $\pi_{\PV}(C)$.  This follows from their Subcurve Inequality which thus arises much earlier in their argument than in Gieseker's. The filtrations used are much the same as in Gieseker's argument (the inequality for a subcurve $Y$ uses essentially the $F_Y$ defined on p.239 of~\cite{Moduli}) but their inequality is first proved in Proposition~5.6 only for subcurves $Y$ no component of which is collapsed by $\pi_{\PV}$.  This inequality, applied to the subcurve $Y$ consisting of all such components, is then used to show that $Y$ must be all of $C$. 

From this point on, the claims and proofs follow closely the pattern of Gieseker's argument with the addition of a check---like that in the last paragraphs of~(\ref{swinarksisetup})---that the marked points must be distinct both from each other and from the nodes. 

The second containment---$J \subset \semi{J}(L)$---corresponds to Corollary~\ref{StabilityAndJss}.\ref{semistablereductionstep}. It's impossible to imitate the semistable replacement argument used to deduce this from the stability of Hilbert points of smooth maps because there are stable maps that cannot be smoothed. In general, $\mgnbardhat$ is reducible, and there are components consisting entirely of nodal maps.  

But, if the set $\semi{J}(L)$ is non-empty, the containment $J \subset \semi{J}(L)$ can be shown by other means.  In~\cite{BaldwinSwinarski}, this is achieved by relying on other constructions of $\mgnbardhat$. Proposition~6 of~\cite{FultonPandharipande} is first used to assert that quotient map $\semi{J} \to \mgnbardhat$ has closed image (cf.~\cite{BaldwinSwinarski}*{Proposition 3.7}). Since this image is necessarily open, it then suffices to know that that $\mgnbardhat$ is connected for which the paper~\cite{KimPandharipande} is invoked. This restricts the results of~\cite{BaldwinSwinarski} to characteristic~$0$.

In \cite{BaldwinPositive}, Baldwin gives a second proof that avoids drawing on these results and works in sufficiently large positive characteristics. An easy argument (her Proposition~2.6) shows that $\semi{J}$ is closed in $J$ by setting up a semi-stable replacement and then invoking the separatedness of the moduli functor of pointed stable maps~(\cite{BehrendManin}*{Proposition~4.1}). The main result, Proposition~2.7, asserting that $J$ is connected is then proved by a strategy similar to that in~\cite{KimPandharipande}. Baldwin begins constructing, for each $1$-ps $\rho$, a family of stable maps that contains the stabilization of the limit under $\rho$ of any suitably generic stable map. It then suffices to show (Proposition~4.5) that there is a connected component of $J$ that contains any such limit map lying in $J$. The approach is similar to that in~\cite{KimPandharipande} but additional work is required because in addition to constructing the necessary family of maps, it is necessary to equip their source curves with embeddings in $Z = \PV \times \PP^{r}$: this is done in Propositions~4.3 and~4.4.

Thus, to complete the construction, it remains to show that there exist stable maps $X$ with semistable Hilbert points and this is done by showing this for maps $X$ with smooth domain $C$. However, instead of checking semistability by a direct calculation, Baldwin uses a diagonal induction. The main step, treated in the next subsection, assumes the stability of $(n-1)$-pointed maps $Y$ of genus $g+1$ with possibly singular domain (in practice, a curve with an elliptic tail) with respect to the corresponding central linearization $L^*$ and deduces the stability of $n$-pointed maps $X$ of genus $g$ with smooth domain, but with respect, not to the central linearization in this new case, but to a perturbation $L$ of it. The induction is then closed by applying the following result, extracted from the first five sections of~\cite{BaldwinSwinarski} (more specifically, from Proposition~2.13, Theorems~3.8 and~5.21 and Corollary~5.22). Once again, I omit the lengthy inequalities defining $\mathbf{L}$.  
\beginS{claim}\label{inductionclosure}
Fix $\nu \ge 10$. For any sufficiently large $m$, there is an open convex set $\mathbf{L}$ of linearizations $(m,\mhat,m')$ containing the central linearization $L^*$ with triple $\bigl(m,\frac{c\nu}{2\nu-1}m,  \frac{\nu}{2\nu-1}m^2\bigr)$ and with the property that  $J \quotient_{L'} \ssl(V) = \mgnbardhat$ for \emph{every} linearization $L'$ with parameters in $\mathbf{L}$ if there exists \emph{any} linearization $L$ with parameters in $\mathbf{L}$ for which the $L$-semistable locus $\semi{J}(L)$ is non-empty.
\end{claim}
The only ingredient in this claim that has not yet been discussed is a variation of \GIT argument that justifies the conclusion for all $L'$ from the non-emptiness of $\semi{J}(L)$ for a single $L$. The set $\mathbf{L}$ is defined by the inequalities that arise as hypotheses in the Potential Stability theorem. For such linearizations, that theorem implies that no map with infinite automorphism group can be Hilbert semistable. Hence, for such $L'$, $\semi{J}(L')= \stable{J}(L')$. But the strictly semistable locus is non-empty for any linearization on a wall. Hence, if $\mathbf{L}$ a convex set of linearizations $L'$ such that the stable and semistable loci with respect to each $L'$ coincide, then $\mathbf{L}$ must lie in a single VGIT chamber. Hence the stable and semistable loci are the same for any 2 linearizations $L$ and $L'$ in $\mathbf{L}$. For more details, see Theorem~3.3.2 of Dolgachev and Hu~\cite{DolgachevHu}, or for a direct and elementary argument~\cite{BaldwinSwinarski}*{Proposition~2.13}.

\subsection{Stability of maps with smooth domains}
\label{stability of smoothmaps}

In view of the preceding discussion, it remains only to verify the Hilbert stability of maps with smooth domains with respect to a suitable linearization $L$ with parameters in $\mathbf{L}$. The  discussion of stability of pointed curves teaches us to expect this to be quite a bit trickier when there are marked points and it is. Baldwin finesses this by the diagonal induction mentioned above and laid out in more detail below.

The base case---stability of maps with smooth, unpointed domains---was treated by Swinarski in~\cite{SwinarskiUnpointed} and is now Theorem~6.5 of~\cite{BaldwinSwinarski}. The basic idea is simple: to see that the arguments used to derive the estimates proving Gieseker's Criterion~\ref{GiesekersCriterion} can be applied to derive the same estimates here. Checking this is straightforward so I'll just sketch the argument briefly. However, this approach works only for $g\ge 2$. This is why the cases $(g,n)$ equal to $(0,0)$ and $(1,0)$ listed at the start of the section cannot be handled. and the lack of the latter as a base eliminates $(g,n)=(0,1)$ from the induction. 

We assume that  $\o_C(1,0)$, the bundle that embeds $C$ in $\PV$, is very ample of degree at least $2g+1$, and that the linear series $\PV$ is complete. All that we can assume about the map $f$ is that it is determined by a basepoint free subseries of the complete linear series of $\o_C(0,1)$ but the triviality of the $\ssl(V)$ action on its sections means that it plays only a secondary role.

Fix a weight filtration $F$ in $V$ and, for a subseries $V'$ in $V$, let $\Vhat'$ be the subseries of $\HO{C}{\o_C(1,1)}$ generated by products of the form $x \cdot y$ where $x\in V'$ and $y \in \HO{C}{\o_C(0,1)}$. Next fix a subsequence  $ 1=j_{0}>j_{1}>\cdots>j_{h}=N$ and large integers $m= (p+1)n$ and $\mhat$. Define $\What^{n}_{k,l}$ by simply replacing each $V$ by the corresponding $\Vhat$ in the discussion leading to~\eqref{giesekerfiltration} and let $\Uhat^n_{k,l}$ be its image in $\HO{C}{\o_C(m,\mhat)}$. 

Now observe first, that the base locus of $\Uhat^1_{k,l}$ equals that of the corresponding subseries $U^1_{k,l} \in \HO{C}{\o_C(p+1,0)}$. This allows us to run the argument for the dimension formula~\eqref{GiesekerDimensionFormula} to see that, for large $n$,  the codimensions of $\Uhat^n_{j,k,l}$ and $U^n_{j,k,l}$ are equal. Next, note that both spaces consist of images of monomials of weight at most $(w_0 + (p-l)w_{j_k}+lw_{j_{k+1}})$, again because $\ssl(V)$ acts trivially on sections of $\o_C(0,1)$. Thus all the ingredients that enter into the estimate for $w_F(m)$ coming from  in~\eqref{giesekerfiltration} are identical for the analogous filtration using the spaces defined by the $\Uhat^n_{j,k,l}$ and we get the same upper bound $w_F(m,\mhat)$.  In view of our hypotheses, this last is negative by Theorem~\ref{StabilityOfSmoothCurves}. 

The real novelty in~\cite{BaldwinSwinarski} is in the ingenious inductive step. Here I'll go into somewhat more detail. We begin by laying out the setup in the notation for stable maps established  above.  The reader will note some similarities to subsection~\ref{elliptictails}.

Fix $X = (C; p_1, \ldots, p_n; f_C:C\to \PP^{r})$, a smooth $n$-pointed stable map of genus $g$ and $Z = (E;q; f_E:E\to \PP^{r})$, a smooth $1$-pointed stable map of genus $1$ chosen so that $f_E$ collapses $E$ to the point $f_C(p_n)$ in $\PP^{r}$. Then let $Y= (D; p_1, \ldots, p_{n-1}; f_D:D\to \PP^{r})$ be the $(n-1)$-pointed stable map of genus $(g+1)$ for which $D=C\cup E$---obtained by gluing $p_n$ on $C$ to $q$ on $E$ to obtain a node $p$ on $D$---and $f_D$ restricts to $f_C$ and $f_E$ respectively on $C$ and $E$. 

Baldwin deduces the stability of the Hilbert point of $X$ with respect to a $1$-ps $\rho_C$ and a suitably chosen linearization $L$ (independent of $\rho_C$) from the stability of the Hilbert point of $Y$ with respect to a $1$-ps $\rho_D$ (non-canonically associated to $\rho_C$) and the central linearization $L^*$. I emphasize that the only element of the \GIT setup that $X$ and $Y$ have in common is the target projective space $\PP^{r}$. They are embedded in different spaces, determine points of different Hilbert schemes with different sets of linearizations and are acted on by different special linear groups. Nonetheless, it's possible to compare the two sides just enough to obtain upper bounds for $L$-weights of bases on the $X$ side from the bounds implied by Hilbert stability for $L^*$-weights of bases on the $Y$ side. That these bounds are sufficient to check the Hilbert stability of $X$ is a minor miracle that depends essentially on the particular choice of linearization $L^*$. The definitions of $\rho_D$ and $L$ will emerge as we review the argument.

First, note that $\o_D(0,1) := f_D^*\bigl(\o_{\PP^{r}}(1)\bigr)$ restricts to 
$\o_C(0,1) := f_C^*\bigl(\o_{\PP^{r}}(1)\bigr)$ on $C$ and to $\o_E(0,1) := f_E^*\bigl(\o_{\PP^{r}}(1)\bigr) =\o_E$  on $E$.
Likewise, for fixed $\nu$ and $c$, let $\o_D(1,0) := \lbpow{\omega_D(\sum_{k=1}^{n-1} p_k)}{\nu}\otimes \o_D(0, c\nu)$, $\o_C(1,0) := \lbpow{\omega_C(\sum_{k=1}^{n} p_k)}{\nu}\otimes \o_C(0, c\nu)$ and $\o_E(1,0) := \lbpow{\omega_E(q)}{\nu}\otimes \o_E(0, c\nu)$. Since $\restrictedto{\omega_D}{C} = \omega_C(p_n)$ and $\restrictedto{\omega_D}{E} = \omega_E(q)$, we find that $\restrictedto{\o_D(1,0)}{C} = \o_C(1,0)$ and
$\restrictedto{\o_D(1,0)}{E} = \o_E(1,0)$. 

Baldwin next translates this in terms of restriction maps. It will be convenient to denote by $\res{S}{T}$, the restriction of sections from $S$ to a subscheme $T$. Set $V_D=\HO{D}{\o_D(1,0)}$ and let $K_C \in V$ and $K_E \in V_D$ be the kernels of the maps $\res{D}{C}(V_D)$ and $\res{D}{E}(V_D)$ given by restricting to $C$ and $E$ respectively. Then,
\beginS{enumerate}
	\item $K_C \cap K_E = \{0\}$.
	\item Any section vanishing at $p$ is uniquely expressible as the sum of elements of $K_C$ and $K_E$.
	\item $\res{D}{C}(K_E)$ may be canonically identified with $\HO{C}{\o_C(1,0)(-p_n)}$.
	\item $\res{D}{E}(K_C)$ may be canonically identified with $\HO{E}{\o_E(1,0)(-q)}$.
\end{enumerate}
We fix, non-canonically, a line $U$ spanned by an element $u \in V_D$ non-zero at $p$ so that $V_D = K_E\oplus U \oplus K_C$ and let $J_C = K_E\oplus U$. By construction, a section $s\in J_C$ is non-zero on $E$ if and only if it is non-zero at $p$.

Next we need to study restrictions of monomials on $Z = \PP(V_D) \times \PP^{r}$. So let's fix a large bidegree $(m,\mhat)$---henceforth we assume all monomials have this bidegree---and write the two pieces of such a monomial as $M\cdot\Mhat$. If $S$ is a set of monomials we will write $\overline{S}$ for its span. Pick any basis $B_D$ of $V_D$ compatible with the decomposition $V_D= J_C \oplus K_C$. Let $\Omega_0$ be the set of monomials $M\cdot\Mhat$ for which $M$ has no factor from $K_C$,  $\Omega_+$ be the set where $M$ has at least one such factor and $\Omega_{++}$ the subset of $\Omega_+$ where no factor from $J_C$ vanishes at $p$. Since $\Omega_0$ and $\Omega_+$ are complementary, we have $\HO{Z}{\o(m,\mhat)} = \overline{\Omega}_0 \oplus \overline{\Omega}_+ $. The heart of Baldwin's argument is:

\beginS{lemma}\label{baldwinclaims}
Restriction to $D$ induces a direct sum decomposition \begin{equation}\label{balddecomp}
\HO{D}{\o_D(m,\mhat)} =\res{Z}{D}(\overline{\Omega}_0) \oplus \res{Z}{D}(\overline{\Omega}_{++})
\end{equation}
and we may identify the first and second summands with the image and kernel, respectively, in the restriction exact sequence
$$ 
0\to \HO{E}{\o_E(m,\mhat)(-q)}\to \HO{D}{\o_D(m,\mhat)} \to \HO{C}{\o_C(m,\mhat)} \to 0\,.
$$
\end{lemma}
\beginS{proof}
Let's first prove a version of~\eqref{balddecomp} with second term $ \res{Z}{D}(\overline{\Omega}_+)$ on the right side. Since $\res{Z}{D}$ is surjective on sections of $\o(m,\mhat)$, such a decomposition follows from the complementarity of $\Omega_0$ and $\Omega_+$ if we check the two terms have trivial intersection. On the other hand, $\res{D}{C}$ is onto and by construction, its kernel contains $\res{Z}{D}(\overline{\Omega}_+)$. So everything will follow if we check that no section $s$ in $\res{Z}{D}(\overline{\Omega}_0)$ can vanish on $C$. 

I claim that section $s$ in $\sym^m(J_C)\otimes\sym^{\mhat}(k^{r+1})$ that vanishes at $p$ also vanishes on $E$. Thus any monomial vanishing on $C$ also vanishes on $D$ as required. Since the claim makes no reference to our basis, we are free to prove it by choosing a more convenient $B$. We want to select coordinates on $J_C$ so that all but one, that we take to be $u$, lie in $K_E$ and on $k^{r+1}$ so that all but one, that we call $\uhat$, vanish at $g(p)$. In such coordinates, the only monomial that does not vanish at $p$ is $u^m \uhat^{\mhat}$. Thus, if $s(P)=0$ then $s$ has no $u^m \uhat^{\mhat}$ term. But then every monomial appearing in $s$ contains a factor from $K_E$ and so vanishes on $E$.

The argument above also shows that we have a decomposition with second term $\res{Z}{D}(\overline{\Omega}_{++})$. Indeed, any monomial in $\Omega_+$ with a factor from $J_C$ vanishing at $p$ vanishes on $E$ and, since it has a factor from $K_C$, it vanishes on $D$.
\end{proof}

We are now ready to define the $1$-ps $\rho_D$ of $\ssl(V_D)$ that Baldwin associates to a $1$-ps $\rho_C$ of $\ssl(V_C)$ where $V_C = \HO{C}{\o_C(1,0)}$. We view $\rho_D$ and $\rho_C$ as weighted bases $B_D$ and $B_C$ as in Section~\ref{numerical} and normalize $\rho_C$ to have total weight $0$. Choose an identification $\sigma: \HO{C}{\o_C(1,0)} \to J_C \subset V_D$ that agrees with the canonical one on the hyperplane $K_E= \HO{C}{\o_C(1,0)(-p_n)}$. The set $B_D$ is defined to be the union of $\sigma(B_C)$ with any basis $B'$ of $K_C$. We assign each element of $\sigma(B_C)$  the weight of the corresponding element of $B_C$ and assign every element of $B'$ the \emph{least} weight of any element of $B_C$ not vanishing at $p_n$ or, equivalently, the least weight $w(n)$ of any element of $\sigma(B_C)$ not lying in $K_E$. For any $k<n$, any coordinate not vanishing at $p_k$ comes from $B_C$ so the least weight $w(k)$ of such a coordinate is the same for both  $B_C$ and $B_D$. Since $\dim(K_C) = \hO{E}{\o_E(1,0)(-q)}= \nu-1$ and we assumed that $\rho_C$ has average weight $0$, the average weight $$\alpha_D = \frac{(\nu-1)w(n)}{\dim(V_D)}\,.$$

The first step is to relate the minimum weights $w_D(m,\mhat)$ and $w_C(m,\mhat)$ of $B_D$- and $B_C$-monomial bases of bidegree $(m,\mhat)$ using Lemma~\ref{baldwinclaims}. The construction of $B_D$ and $B_C$ means that a $B_C$-basis of $\HO{C}{\o_C(m,\mhat)}$ has the same weight as the part of a $B_D$-basis of $\HO{D}{\o_D(m,\mhat)}$ spanning the first term $\res{Z}{D}(\overline{\Omega}_0)$ of~\eqref{balddecomp}. On the other hand, all the factors of any monomial in the second term of~\eqref{balddecomp} either come from $K_C$ or from elements of $B_C$ not vanishing at $p_n$ and hence have weight at least $w(n)$. Putting these remarks together:	
$$ w_D(m, \mhat) \ge w_C(m,\mhat) + \Bigl(\hO{D}{\o_D(m,\mhat)} - \hO{C}{\o_C(m,\mhat)}\Bigr)m w(n)$$

On the other hand, since we are inductively assuming that $Y$ is Hilbert stable, Proposition~\ref{NumericalCriterionForStableMaps} implies that
$$ w_D(m, \mhat) +  m'\sum_{k=1}^{n-1}w(k) <  \bigl(m\p_D(m,\mhat) + m'(n-1)\bigr)\alpha_D\,. $$

Combining these we get an upper bound for $w_C(m,\mhat)$ not involving $ w_D(m, \mhat)$. To unwind this, we first use Riemann-Roch to evaluate
$$P_D(m,\mhat) = \hO{D}{\o_D(m,\mhat)} = 
			m\nu\bigl(2g+(n-1)+c\dhat\bigr)+
			\mhat\dhat-g$$
---which also computes $\dim(V_D)= \hO{D}{\o_D(1,0)}=P_D(1,0)$---and 
$$\hO{C}{\o_C(m,\mhat)} = 
			m\nu\bigl(2g-2)+n+c\dhat\bigr)+
				\mhat\dhat-g+1\,.$$
Then we use the hypothesis that $Y$ is Hilbert stable with respect to the central linearization $L^*$ for which $\frac{\mhat}{m} = \frac{c\nu}{2\nu-1}$ and $\frac{m'}{m^2} = \frac{\nu}{2\nu-1}$ to eliminate $m'$ and $\mhat$. Baldwin next makes all these substitutions, and carries out some lengthy, but completely elementary, algebraic simplifications that I will omit, noting only that the choice of $L^*$ leads to a miraculous series of cancellations. When the dust settles, we obtain,
$$ w_C(m,\mhat) + m^2 \frac{\nu}{2\nu-1} \sum_{k=1}^n w(k)-m\gamma w(n) \le 0$$
with
$$\gamma := 1-\frac{(\nu-1)g}{(2\nu-1)g+\nu(n-1+c\dhat\,)}\,.$$
In particular, $\frac{1}{2} \le \gamma \le 1$.

Repeating this argument with each $p_k$ in the distinguished role played by $p_n$ above, and then averaging gives, 
$$
w_C(m,\mhat) + m''\sum_{k=1}^{n}w(k)\le 0
\text{\quad with \quad} m''=\frac{\nu m^2}{2\nu-1}-\frac{\gamma m}{n}\, .$$

Applying Proposition~\ref{NumericalCriterionForStableMaps} again, the upshot is that if the $(n-1)$-pointed stable map $Y$ of genus $(g+1)$ is Hilbert stable for the central linearization $L^*$ with parameters $(m,\mhat, m')$, then the $n$-pointed stable map $X$ of genus $g$ is Hilbert stable for the linearization $L$ with parameters $(m,\mhat, m'')$. 

Because $|\gamma | < 1$, the difference $|\frac{m''}{m^2} - \frac{\nu}{2\nu-1}|$ is at most $\frac{1}{m}$ and it is straightforward to check that $L$ again lies inside the set $\mathbf{L}$ of linearizations in Claim~\ref{inductionclosure}: for details, see~\cite{BaldwinSwinarski}*{pp.73--74}. This completes the inductive step.


\section{Open problems}
\label{open}

My goal here is to list some questions about Hilbert stability of curves whose answers I would be interested in knowing. Most complement or refine constructions discussed in the previous sections.

\subsection{Combining the Span Lemma with convexity arguments} \label{improvedspan}
The geometric ingredients that go into the proof of stability of smooth curves---estimates from Riemann-Roch and Clifford's Theorem---are the same for unpointed curves. Yet the proof of the former here  takes just over 1 page and the proof of the latter in~\cite{SwinarskiThesis} takes almost 25. The difference is due to the lack of an analogue, in the pointed case, of the convexity argument in Lemma~\ref{combinatorialclaim}, coming from the fact that what is being sought can be viewed as the optimum of a linear program. This point of view is pushed even further in \cite{MorrisonRuled}*{Theorem~4.1}.

A review of Swinarski's argument reveals that it involves estimating areas under piecewise linear and step functions but while the choices are made to achieve local convexity with respect to each base point, they preclude adjusting the weights to obtain any global convexity. 

Is there a way of using the Span Lemma~\ref{spanlemma} to get an estimate for $e_F$ that is a convex function of the weights (or even a sum of such functions)? My hunch is that the answer is yes and that such an estimate would lead to a much simpler proof of the stability of Hilbert points of smooth pointed curve. I think such a result would have other applications: see (\ref{directnodal}). 

\subsection{Pointed pseudostability and other variants} 
\label{pointedvaiants}
In running the log minimal model program for $\mgnbar$, you'd expect to encounter most of the variant moduli problems discussed in Section~\ref{schubert} that arise in running it for $\mgbar$ such as Schubert's pseudostable curves and the $c$- and $h$-semistable curves (having tacnodes as well as cusps) of Hassett and Hyeon, and you'd expect coarse moduli spaces for these problems to arise as GIT quotients of $\nu$-canonical loci for small $\nu$ as for unpointed curves. What are the small-$\nu$ quotients that arise?

Such questions already arise in the log minimal model program for $\mgbar$.  I'll sketch just one example from a study in progress, by Hyeon and Lee, of the map $\mgbar(\frac{7}{10}-\epsilon\bigr)\to \mgbar(\frac{2}{3}\bigr)$---see also~\ref{smallm}. They need to address questions about bicanonical stability of $1$-pointed curves because to understand the restrictions of the relevant log canonical classes to the stratum $\Delta_2$ in the boundary of the moduli space of h-semistable curves discussed in~\ref{logminimal}, it is necessary to understand their restrictions to the pointed factors $\mijbar{2}{1}$ and $\mijbar{g-2}{1}$. To do so, they consider curves $C$ in $\mijbar{2}{1}$ embedded by $\lbpow{\omega_C}{2}(2p)$ with $p$ the marked point. Here their calculations show that if $p$ lies on an elliptic bridge $E$ (a genus $1$ component meeting a genus $0$ component in $2$ nodes), then $C$ is Hilbert unstable and in the quotient is replaced by a curve in which $E$ is contracted to a tacnode.

I also want to mention David Smyth's paper~\cite{SmythThesis} which appeared while this article was in revision, and which touches on a number of ideas discussed here although its does not use \GIT quotients. In it, he constructs a proper Deligne-Mumford moduli stack $\Mijbar{1}{n}(m)$ for $m$-stable $n$-pointed curves of genus $1$ whenever $m<n$. His construction does not employ \GIT so I will not even precis it here, but the $m$-stability condition deserves a brief mention. Smyth first defines an elliptic $k$-fold point to be one locally isomorphic to an ordinary cusp for $k=1$, a tacnode for $k=2$ and the intersection of $k$ lines in $\AA^{k-1}$ for $k\ge 3$.  Then $m$-stability for $(C, p_1, \ldots, p_n)$ permits only nodes and elliptic $k$-fold points for $k \le m$ as singularities, and requires first that the sum of the number of marked points on a connected subcurve of genus $1$ and the number of points in which it meets the rest of the curve be at least $m+1$, and second, a vanishing condition---that $\HO{C}{\omega_C^{\vee}(-\sum_{i=1}^n p_i)}=\{0\}$---that I won't go into here. Thus, were the case $n=0$ to make sense in genus $1$, the first two conditions for $m=1$ and $m=2$ respectively would match the notions of pseudostability and $h$-stability in Section~\ref{schubert}.

Two other recent papers in this area that make contact with the ideas here but do not use \GIT qoutients trhat 'd like to flag for the interested reader are Dawei Chen's study~\cite{ChenLogMinimal} of $\overline{\scr M}\sb {0,0}(\mathbb P^{3}, 3)$ and Matt Simpson's work~\cite{SimpsonLogCanonical} on log canonical models of $\mijbar{0}{n}$.

\subsection{Direct proof of Hilbert stability for nodal stable curves}\label{directnodal}
This problem is in the nature of an embarrassing lacuna in the subject. We know which nodal stable curves have stable Hilbert points. The Potential Stability Theorem~\ref{PotentialStabilityTheorem} gives necessary conditions, at least when the degree is large enough compared to the genus, and Caporaso~\cite{CaporasoCUP} provides a fairly complete converse. Her arguments are, like Gieseker's via semi-stable replacement. The annoying fact is that no one has ever verified the Numerical Criterion by estimating weights of bases.

As we remarked above, Gieseker's Criterion~\ref{giesekerscriterion} is too weak. My hunch is that the Span Lemma~\ref{spanlemma} might provide the extra tool needed to find sufficiently sharp estimates.

We can get a feel for the difference going back to Swinarski's Example~1---see~\eqref{exampleone}---and supposing that the points $p_i$ are nodes. Then the estimates given there continue to apply with the difference that all the codegrees $e_i$ are doubled. Thus the estimate $\epsilon_F =3$ for $e_F$ in Gieseker's Criterion is $3$, while the right hand side is slightly larger than $2$. But once again, estimating the subspace of $\HO{C}{\lbpow{L}{m}}$ weight at most $\frac{2m}{6}$ not by $V_1^m$ as in Gieseker's criterion but by the span of $V_1^m$ and $\bigl(V_0 V_2\bigr)^{m/2}$ reduces the estimated codimension of this space from $2m$ to $m$ and the estimate for $e_F$ by $\frac{2}{3}$. In fact, just as applying the Span Lemma systematically reduced the estimate for $e_f$ in Swinarksi's original example from $\frac{3}{2}$ to $1$, so in this modified example we get a reduction from $3$ to $2$.

Working this out in general would almost certainly be more delicate than handling pointed curves, so it probably makes sense to try to squeeze more out of the Span Lemma as proposed in~\ref{improvedspan} first. 

\subsection{Stability of canonical models of stable curves} \label{canonical}

In Lemma~\ref{StabilityOfCanonicalNonHECurves}, we proved the asymptotic Hilbert stability of canonical models of smooth non-hyperelliptic curves. There are several interesting questions about the stability of canonical models not covered by this.

The most obvious concerns rational ribbons in the sense of Bayer and Eisenbud~\cite{BayerEisenbudRibbons} that arise as the flat limits of canonical models of non-hyperelliptic curves as they approach the hyperelliptic locus. What can we say about their Hilbert points? If they are semistable, it would be possible to apply the Cornalba-Harris Theorem~\ref{cornalbaharristheorem} in a unified way to any generically smooth family of curves. Stoppino~\cite{StoppinoSlopeInequalities} shows that equality~\eqref{smoothnonheinequality} is sharp exactly for families with generic fiber hyperelliptic. In the same spirit, Moriwaki~\cite{MoriwakiRelative}*{Theorem~5.1} gives inequalities necessary and sufficiently for a divisor on $\mgbar$ to meet effectively all curves not lying in $\Delta$ and the curves that show the necessity all lie in the closure of the hyperelliptic locus. Both these results suggest that such ribbons are Hilbert strictly semistable. 

More generally, the work of Hassett and Hyeon discussed in (\ref{logminimal}) indicates that canonical quotients will arise as log minimal models $\mgbar(\alpha)$ for lower values of $\alpha$. Do the Chow and Hilbert quotients have a modular interpretation in terms of some class of curves like the pseudostable curves of~\cite{Schubert} or the $c$- and $h$-semistable curves of~\cite{HassettHyeonFlip}?

\subsection{Checking $m$-Hilbert stability for small values of $m$} \label{smallm}

All the methods for verifying Hilbert stability that we use here are fundamentally asymptotic. That is, they show $m$-Hilbert stability for all sufficiently large $m$ by showing that the leading coefficient $e_{F}$ of Lemma~\ref{AsymptoticNumericalCriterion} is bounded away from $0$ from below. The techniques used to get these estimates for a general  filtration $F$ or $1$-ps $\rho$ seem to require the freedom both to assume that $m$ is large both to invoke vanishing theorems and use Riemann-Roch to compute dimensions and to ignore various terms of order $m$.  It often happens that we can verify that a curve is $m$-Hilbert unstable for a small $m$, since this requires us only to compute weights for a single $1$-ps $\rho$ which is often possible either by an exact deductive calculation or by a symbolic computation. But to check semistability, we must handle all $\rho$. The formula~\eqref{mtwothree} of Hassett and Hyeon is a good example of a criterion that can be used in both these ways to check instability but that seems difficult to apply to check semistability.

So an interesting, if apparently difficult, problem is to find methods for checking $m$-Hilbert semistability for small $m$. Any such methods would find an immediate application in the log minimal model program. Unpublished calculations of Hassett using the work of Gibney, Keel and the author~\cite{GibneyKeelMorrison} on the  $F$-conjecture lead to predictions of critical values of $\alpha$ at which the birational type of $\mgbar(\alpha)$ will change. 

I conclude with two examples. These methods predict that, for $g \ge 4$, the next critical value below $\frac{7}{10}$ is $\alpha=\frac{2}{3}$ when $K_{\Mgbar}+\frac{2}{3}\delta= \frac{1}{3}(39\lambda-4\delta)$ has slope $9.75$. Comparing this with the polarization formula~\eqref{polarizationformula}, we see that a polarization with exactly this slope is predicted for the quotient of the degree $m=6$ Hilbert scheme of bicanonically embedded curves. In work in progress, Hyeon and Lee use results of Rulla~\cite{RullaThesis} to understand the pullback of $K_{\Mgbar}+\frac{2}{3}\delta$ to $\mijbar{2}{1}$, Hyeon and Lee conclude that this class should contract the locus of genus $2$ Weierstrass tails (curves having a genus two subcurve meeting the rest of the curve in a node that is a Weierstrass point of the tail). 

Hassett and Hyeon have shown that such a curve $C$ is $6$-Hilbert
strictly semistable with respect to a $\rho$ analogous to that treated in~\ref{elliptictails} with weights given by the ramification sequence at the point of attachment. Further, they show that the flat limit is obtained by replacing the Weierstrass tail by a rational tail carrying a ramphoid cusp singularity (one analytically isomorphic to $y^2=x^5$). The same calculation shows that curves with a Weierstrass tail are $5$-Hilbert unstable. This suggests that there is flip
\setdiagram{flipsixdiagram}{nohug,grid=flipgrid}{
Q^{\infty}& & & & Q^{5}\\
& \rdTo^{\Psi}& &\ldTo^{\Psi^+} & \\
& & Q^{6}& & \\
}%
with $Q^m$ the $\ssl$-quotient of the Hilbert scheme of $2$-canonically embedded curves linearized by $\Lambda_m$ and  $Q^{\infty}$ the asymptotic limit (given by taking an $m \gg 0$).

Similar considerations predict that the next critical value is  $\alpha=\frac{19}{29}$ when $K_{\Mgbar}+\frac{19}{29}\delta= 13(\lambda-\frac{3}{29}\delta)$ has slope $\frac{29}{3}$ which, by~\eqref{polarizationformula}, occurs for 
the quotient of the degree $m=\frac{9}{2}$ Hilbert scheme of bicanonically embedded curves. Of course, this Hilbert scheme makes no sense but, by clearing denominators, the polarization $\Lambda_{\frac{9}{2}}$ does---or, one can simply work with $m=4$. Here, Hassett has made calculations that show that curves in $\Delta_2$ are $m=\frac{9}{2}$ strictly-semistable and $m=4$ unstable with respect to the $\rho$ analogous to that in~\ref{elliptictails} having weights given by the ramification sequence at the point of attachment on the genus $2$ tail (no longer, in general, a Weierstrass point). This suggests that here $\Delta_2$ gets contracted.


\linespread{.965}\normalfont\selectfont
\ifbells
	
\fi

\end{document}